\documentclass[12pt]{article}
\usepackage{amsmath,amsthm,amssymb}

\newtheorem{theo}{Theorem}[section]
\newtheorem{lemma}[theo]{Lemma}
\newtheorem{coro}[theo]{Corollary}
\newtheorem{pro}[theo]{Proposition}

\numberwithin{equation}{section}

\theoremstyle{definition}

\newtheorem*{example*}{Example}

\newtheorem*{remark*}{Remark}

\newcommand{\bZ}{{\mathbb Z}}
\newcommand{\bN}{{\mathbb N}}
 \DeclareMathOperator{\End}{End}

 \DeclareMathOperator{\Jalg}{Jalg}
 \DeclareMathOperator{\Jideal}{Jideal}
 \DeclareMathOperator{\ideal}{ideal}
 \DeclareMathOperator{\trace}{trace}

\newcommand{\frg}{{\mathfrak g}}

\newcommand{\frh}{{\mathfrak h}}

\newcommand{\subo}{_{\bar 0}}
\newcommand{\subuno}{_{\bar 1}}

\title{Maximal subalgebras of Jordan superalgebras.}

\author{Alberto Elduque\thanks{The first and second authors
have been supported by the Spanish Ministerio de Educaci\'on y
Ciencia and FEDER (MTM 2004-081159-C04-02), and the second and third
authors by the Comunidad Aut\'onoma de La Rioja. The first author
also acknowledges support by the Diputaci\'on General de Arag\'on
(Grupo de Investigaci\'on de
\'Algebra).}
\\
{\small Departamento de Matem\'aticas}\\
{\small Universidad de Zaragoza}\\
{\small 50009, Zaragoza. Spain}
\and
Jes\'us Laliena{${}^*$} and Sara Sacrist\'an{${}^*$}
\\{\small Departamento de Matem\'aticas y Computaci\'on}\\
{\small Universidad de La Rioja}\\
{\small 26004, Logro\~no. Spain}}

\date{\quad}

\begin{document}
\maketitle\vspace{-1.5cm}

\begin{abstract}
The maximal subalgebras of the finite dimensional simple special
Jordan superalgebras over an algebraically closed field of
characteristic $0$ are studied.  This is a continuation of a
previous paper by the same authors about maximal subalgebras of
simple associative superalgebras, which is instrumental here.
\end{abstract}


\bigskip

\section{Introduction.}

\bigskip

Finite dimensional simple Jordan superalgebras over an algebraically
closed field of characteristic zero were classified by V.~Kac in
1977 \cite{Kac}, with one missing case that was later described by
 I.~Kantor in 1990 \cite{Kan}. More recently
M.~Racine and E.~Zelmanov \cite{RaZe} gave a classification of
finite dimensional simple Jordan superalgebras over arbitrary fields
of characteristic different from 2 whose even part is semisimple.
Later, in 2002, C.~Mart\'{\i}nez and E.~Zelmanov \cite{Ma-Ze}
completed the remaining cases, where the even part is not
semisimple.

Here we are interested in describing the maximal subalgebras of the
finite dimensional simple special Jordan superalgebras with
semisimple even part over an algebraically
 closed field of characteristic zero. Precedents of this work
are the papers of E.~Dynkin in 1952  (see \cite{Dy1}, \cite{Dy2}),
where the maximal subgroups of some classical groups and  the
maximal subalgebras of semisimple Lie algebras are classified, the
papers of M.~Racine (see \cite{Ra1}, \cite{Ra2}), who classifies the
maximal subalgebras of finite dimensional central simple algebras
belonging to one of the following classes: associative, associative
with involution, alternative and special and exceptional Jordan
algebras; and the paper
 by the first author in 1986
(see \cite{El}), solving the same question for central simple Malcev
algebras.

In a previous work \cite{El-La-Sa}, the authors described the
maximal subalgebras of finite dimensional central simple
superalgebras which are either associative or associative with
superinvolution. The results obtained there will be useful in the
sequel. The maximal subalgebras of the ten dimensional Kac Jordan
superalgebra are determined in \cite{El-La-SaKAC}.

\smallskip

 First of all, let us recall some basic facts.
A \emph{superalgebra} over a field $F$ is just a $\bZ_2$-graded
algebra $A=A\subo\oplus A\subuno$ over $F$ (so $A_\alpha
A_\beta\subseteq A_{\alpha+\beta}$ for $\alpha,\beta\in\bZ_2$). An
element $a$ in $A_\alpha$ ($\alpha=\bar 0,\bar 1$)  is said to be
\emph{homogeneous} of degree $\alpha$ and the notation $\bar
a=\alpha$ is used. A superalgebra is said to be \emph{nontrivial} if
$A\subuno\ne 0$ and \emph{simple} if $A^2\ne 0$ and $A$ contains no
proper graded ideal.

An \emph{associative superalgebra} is just a superalgebra that is
associative as an ordinary algebra. Here are some important
examples:

\bigskip

\begin{enumerate}

\item[{\rm a)}]  $A=M_n(F)$, the algebra of $n\times n$ matrices over $F$, where

$A_{\bar{0}} =\left\{ \begin{pmatrix} a & 0 \\ 0 & b \end{pmatrix}
: a\in M_r(F), b\in M_s(F)\right\},$

$ A_{\bar{1}}= \left\{ \begin{pmatrix}  0 & c\\ d & 0\end{pmatrix}
: c\in M_{r\times s}(F), d\in M_{s\times r}(F)\right\},$

with $r+s =n.$ This superalgebra is denoted by $M_{r,s}(F).$

\bigskip

\item[{\rm b)}] The subalgebra $A=A_{\bar{0}}\oplus A_{\bar{1}}$ of $M_{n,n}(F)$, with

$A_{\bar{0}}= \left \{ \begin{pmatrix} a & 0\\ 0 & a \end{pmatrix}
: a\in M_n(F)\right\}, \quad A_{\bar{1}} = \left\{ \begin{pmatrix} 0 & b\\
b & 0
\end{pmatrix} : b\in M_n(F)\right \}.$

This superalgebra is denoted by $Q_n(F).$

\bigskip

Over an algebraically closed field, these two previous examples
exhaust the simple finite dimensional associative superalgebras, up
to isomorphism.

\bigskip

 \item[{\rm c)}] The \emph{Grassmann superalgebra}:
\[
G= alg \langle 1,e_1, e_2, \dots \ : e_i^2=0= e_ie_j +e_je_i  \
\forall i,j=1,2,\dots \rangle
\]
 over a field $F$, with the grading $G=G_{\bar{0}}\oplus
G_{\bar{1}},$ where $G_{\bar{0}}$ is the vector space spanned by the
products of an even number of $e_i$'s, while $G\subuno$ is the
vector subspace spanned by the products of an odd number of $e_i$'s.
(The product of zero $e_i$'s is, by convention, equal to $1$.)

\end{enumerate}

\bigskip

Following standard conventions, given a superalgebra $A=A\subo\oplus
A\subuno$, the tensor product $G\otimes A$, where $G$ is the
Grassmann superalgebra, becomes a superalgebra with the product
given by $(g\otimes a)(h\otimes b)=(-1)^{\bar a\bar h}gh\otimes ab$
for homogeneous elements $g,h\in G$ and $a,b\in A$, and grading
given by $(G\otimes A)\subo=G\subo\otimes A\subo\,\oplus\,
G\subuno\otimes A\subuno$, $(G\otimes A)\subuno=G\subo\otimes
A\subuno\,\oplus\, G\subuno\otimes A\subo$. Its even part
$G(A)=(G\otimes A)\subo$ is called the \emph{Grassmann envelope} of
the superalgebra $A$. Moreover, the superalgebra $A$ is said to be a
superalgebra in a fixed variety if $G(A)$ is an ordinary algebra
(over $G\subo$) in this variety. In particular, $A$ is a Jordan
superalgebra if and only if $G(A)$ is a Jordan algebra.

It then follows that over fields of characteristic $\ne 2,3$,  a
superalgebra $J=J\subo\oplus J\subuno$ is a Jordan superalgebra if
and only if for any homogeneous elements $a,b,c$ in $J$:
\[
L_a b= (-1)^{\bar a\bar b} L_b a,
\]
where $L_a$ denotes the multiplication by $a$, and
\begin{equation}
\begin{split}
L_{a}L_{b}L_{c}&
  + (-1)^{\bar a\bar b + \bar a\bar c + \bar b\bar c}L_{c}L_{b}L_{a}
   + (-1)^{\bar b\bar c }L_{(ac)b}  \\
 &=L_{ab}L_{c}+ (-1)^{\bar b\bar c}L_{ac} L_{b}+
    (-1)^{\bar a\bar b+ \bar a\bar c}L_{bc}L_{a}  \\
 &=(-1)^{\bar a\bar b}L_{b}L_{a} L_{c}+
   (-1)^{\bar a\bar c+ \bar b\bar c}L_{c}L_{a}L_{b}
      + L_{a(bc)}  \\
 &=(-1)^{\bar a\bar c+ \bar b\bar c}L_{c}L_{ab}
    + (-1)^{\bar a\bar b}L_{b} L_{ac}+
    L_{a}L_{bc} .
\end{split}
\end{equation}

\bigskip

Let $A$ be a superalgebra. A \emph{superinvolution} is  a graded
linear map $* \colon A \to A$ such that $x^{**} =x$, and $(xy)^*=
(-1)^{\bar x\bar y}y^*x^*$, for any homogeneous elements $x,y$ in
$A$.

The simplest examples of Jordan superalgebras over a field of
characteristic $\ne 2$ are the following:

\begin{enumerate}
\item[{\rm i)}] Let $A=A_{\bar{0}}+A_{\bar{1}}$ be an associative
superalgebra. Replace the associative product in $A$ with the new
one: $x\circ y = {\frac{1} {2}}(xy + (-1)^{\bar x\bar y}yx)$. With
this product $A$ becomes a  Jordan superalgebra, denoted by $A^+$.

\item[{\rm ii)}] Let $A$ be an associative superalgebra with superinvolution
$*$. Then the subspace of hermitian elements $H(A,*) = \{ a \in A :
a^* =a\}$ is a subalgebra of $A^{+}$.

\end{enumerate}

\bigskip

In fact, if a Jordan superalgebra $J$  is a subalgebra of $A^+$ for
an associative superalgebra $A$,  $J$ is said to be \emph{special}.
Otherwise $J$ is said to be \emph{exceptional}. Any graded Jordan
homomorphism $\sigma \colon J \to A^+$ is called a
\emph{specialization}. So $J$ is special if there exists a faithful
specialization of $J$. Otherwise, $J$ is exceptional. Both examples
i) and ii) given above are examples of special Jordan superalgebras.

A specialization $u \colon J \to U^+$ into an associative
superalgebra $U$ is said to be
 \emph{universal} if
the subalgebra of $U$ generated by $u(J)$ is $U$, and for any
arbitrary specialization $\varphi \colon J \to A^+$, there exists a
homomorphism of associative superalgebras $\chi \colon U \to A$ such
that $\varphi =\chi \circ u$. The superalgebra $U$ is called the
\emph{universal enveloping algebra} of $J$.

\bigskip

\emph{In the sequel only  finite dimensional Jordan superalgebras
over an algebraically closed field of characteristic zero will be
considered.}

\bigskip

We recall the classification of the nontrivial simple Jordan
superalgebras given by V.~Kac \cite{Kac} and completed by I.~Kantor
\cite{Kan}.

\bigskip

\begin{enumerate}

\item[{\rm 1)}] $J=K_3,$ the \emph{Kaplansky superalgebra}:

 $J_{\bar{0}}=Fe,\quad J_{\bar{1}}=Fx+Fy,$
 \quad $e^2=e, \quad e \cdot x=\frac{1}{2} x,
 \quad  e \cdot y = \frac{1}{2}y, \quad x \cdot y=e$.

\bigskip

\item[{\rm 2)}]  The one-parameter family of superalgebras $J=D_t,$
with $t\in F \setminus\{0\}$:

 $J_{\bar{0}}=Fe +Ff,$\quad  $J_{\bar{1}} = Fu+Fv$

 $e^2 =e, \quad f^2=f, \quad  e \cdot f=0,\quad
 e \cdot u=\frac{1}{2}u, \quad e \cdot v=\frac{1}{2}v, \quad
f \cdot u=\frac{1}{2}u,$

$ f \cdot v=\frac{1}{2}v, \quad u \cdot v=e+tf.$

Note that $D_t\cong D_{1/t}$, for any $t\neq 0.$

\bigskip

\item[{\rm 3)}] $J=K_{10},$ the \emph{Kac superalgebra}. This is a
ten dimensional Jordan superalgebra with six dimensional even part.
(See \cite{Hogben-Kac}, \cite{Ki}, \cite{N-B-E} or
\cite{El-La-SaKAC} for details).

\bigskip

\item [{\rm 4)}] Let $V= V_{\bar{0}}\oplus V_{\bar{1}}$ be a graded vector space over $F,$ and let
$(\ , \ )$ be a nondegenerate supersymmetric bilinear superform on
$V,$ that is, a nondegenerate bilinear map which is symmetric on
$V_{\bar{0}},$ skewsymmetric on $V_{\bar{1}},$ and $V_{\bar{0}}$ and
$V_{\bar{1}}$ are orthogonal relative to $(\ , \ ).$ Now consider
$J_{\bar{0}}= Fe+ V_{\bar{0}},$ $J_{\bar{1}}=V_{\bar{1}}$ with $e
\cdot x=x,$ $v \cdot w=(v,w)e,$ for any $x\in J$ and $v,w\in V$.
This superalgebra $J$ is called the \emph{superalgebra of a
superform}. If $\dim V_{\bar{0}} =1$ and $\dim V_{\bar{1}}=2,$ the
superalgebra of a superform is isomorphic to $D_t$ with $t=1.$

\bigskip

\item [{\rm 5)}] $A^{+},$ with $A$ a finite dimensional simple associative
superalgebra, that is, either $A=M_{r,s}(F)$ or $A=Q_n(F)$. Note
that $M_{1,1}(F)^+$ is isomorphic to $D_{-1}$.

\bigskip

\item [{\rm 6)}] $H(A,*),$ where $A$ and $*$ are of one of the following types:

{\parindent =1em   i) $A=M_{n,n}(F),$ $*\colon \begin{pmatrix} a &
b \\ c & d\end{pmatrix} \to  \begin{pmatrix} d^t & -b^t \\ c^t &
a^t \end{pmatrix}.$

 ii) $A=M_{n,2m}(F),$ $* \colon \begin{pmatrix} a & b \\
c & d \end{pmatrix} \to \begin{pmatrix} a^t & c^tq \\
-q^tb^t & q^td^tq \end{pmatrix},$ where $q =\begin{pmatrix} 0 &
I_m
\\ -I_m & 0 \end{pmatrix}.$}

The first one is called the {\sl transpose superinvolution} and
$H(A,*)$ is denoted then by $p(n),$ and the second one the {\sl
orthosymplectic superinvolution} and $H(A,*)$ is denoted in this
case by $osp_{n,2m}$. The isomorphisms $D_{-2} \cong
 D_{-1/2} \cong osp_{1,2}$ are easy to prove.

\bigskip

\item [{\rm 7)}] Let $G$ be the Grassmann superalgebra. Consider the following
product in
$G$:
\[
\{f,g\} = \sum_{i=1}^n (-1)^{\bar f} \frac{\partial f}{\partial
e_i}\frac{\partial g}{\partial e_i},
\]
and build the vector space, sum of two copies of $G$: $J=G+Gx$, with
 the product in $J$ given by
\[
a(bx) = (ab)x, \quad (bx)a=(-1)^{\bar a}(ba)x, \quad
(ax)(bx)=(-1)^{\bar b}\{a,b\}.
\]

\noindent Finally  take the following grading in $J$:
$J_{\bar{0}}=G_{\bar{0}}+G_{\bar{1}}x,
J_{\bar{1}}=G_{\bar{1}}+G_{\bar{0}}x.$ This superalgebra is called
the {\it Kantor double of the Grassmann algebra} or the \emph{Kantor
superalgebra}.

\end{enumerate}

The 10-dimensional Kac superalgebra and the Kantor superalgebra are
the unique exceptional superalgebras in the above list (see \cite
{McC2} and \cite {Sht}). Note that the Kaplansky superalgebra is the
unique nonunital simple superalgebra.

\bigskip

 Let $J$ be a non unital Jordan superalgebra, the unital hull of $J$ is defined to be $H_F(J)=
J + F\cdot 1$, where $1$ is the formal identity and $J$ is an ideal
inside $H_F(J)$. In \cite {Ze} E. Zelmanov determined a
 classification theorem for finite dimensional  semisimple Jordan superalgebras.

\begin{theo}\label{th:semisimple} \emph{(E. Zelmanov)}

Let $J$ be a finite dimensional Jordan superalgebra over a field $F$
of characteristic not 2. Then $J$ is semisimple if and only if $J$
is a direct sum of simple Jordan superalgebras and unital hulls
$H_K(J_1 \oplus \dots \oplus J_r)= (J_1 \oplus \dots \oplus J_r) + K
\cdot 1$ where $J_i$ are non unital simple Jordan superalgebras over
an extension $K$ of $F$.

\end{theo}

\bigskip

The maximal subalgebras of the Kac Jordan superalgebra (type 3)
above) have been determined in \cite{El-La-SaKAC}.  Our purpose in
this paper is to describe the maximal subalgebras of the simple
special Jordan superalgebras (types 1), 2), 4), 5) and 6)). This is
achieved completely for the simple Jordan superalgebras of types 1),
2) and 4). For types 5) and 6)  the results are not complete and
some questions arise.

\bigskip

In what follows the word subalgebra will always be used in the
graded sense, so any subalgebra is graded.

\bigskip

First note that any maximal subalgebra $B$ in a simple unital Jordan
superalgebra $J$, with identity element $1$, contains the identity
element. Indeed, if $1\notin B$, the algebra generated by $B$ and
$1$: $B+F \cdot 1$, is the whole $J$ by maximality. So $B$ is a
nonzero graded ideal of $J$, a contradiction with  $J$ being simple.
Therefore  $1\in B$.

\bigskip

The paper is organized as follows. Section 2 deals with the easy
problem of determining the maximal subalgebras of the Kaplansky
superalgebra, the superalgebras $D_t$ and the Jordan superalgebras
of superforms. Then Section 3 will collect some known results on
universal enveloping algebras and will put them in a way suitable
for our purposes. Sections 4 and 5 will be devoted, respectively, to
the description of the maximal subalgebras of the simple Jordan
superalgebras $A^+$ and $H(A,*)$, for a simple finite dimensional
associative algebra $A$, and a superinvolution $*$.

\section {The easy cases.}

\bigskip

Let us first describe the maximal subalgebras of the simple Jordan
superalgebras of types 1), 2), and 4) in section 1. The proof is
straightforward.

\begin {theo}\label{th:easy}

\begin{enumerate}

\item[{\rm (i)}] Let $J=K_3$ be the Kaplansky superalgebra. A subalgebra $M$ of
$J$ is maximal if and only if $M=J_{\bar{0}}\oplus M_{\bar{1}}$
where $M_{\bar{1}}$ is a vector subspace of $J_{\bar{1}}$ with $\dim
M_{\bar{1}}=1.$

\item[{\rm (ii)}] Let $J=D_t$ with $t\neq 0.$ A subalgebra
$M$ of $J$ is maximal if and only if either $M= J_{\bar{0}} \oplus
M_{\bar{1}}$ where  $M_{\bar{1}}$ is a vector subspace of
$J_{\bar{1}}$ with $\dim M_{\bar{1}} =1,$ or if $t=1,$ $M=F \cdot
1+J_{\bar{1}}.$

\item[{\rm (iii)}] Let $J$ be the Jordan superalgebra of a
nondegenerate bilinear superform. A subalgebra $M$ of $J$ is maximal
if and only if either $M= J_{\bar{0}} \oplus M_{\bar{1}}$ where
$M_{\bar{1}}$ is a vector subspace and $\dim M_{\bar{1}}= \dim
J_{\bar{1}}-1,$ or $M= (F \cdot 1+M_{\bar{0}})\oplus J_{\bar{1}}$
where $M_{\bar{0}}$ is a vector subspace and $\dim M_{\bar{0}}= \dim
V_{\bar{0}}-1.$

\end{enumerate}

\end{theo}

\bigskip

Note that item (ii) in Theorem \ref{th:easy} above cover the maximal
subalgebras of $M_{1,1}(F)^+\cong D_{-1}$ and of $osp_{1,2}\cong
D_2$.

\bigskip

\section{Universal enveloping algebras.}

In order to determine the maximal subalgebras of the remaining
simple special Jordan superalgebras, some previous results are
needed.

Given an associative superalgebra $A$ and a subalgebra $B$ of the
Jordan superalgebra $A^+$, $B'$ will denote the (associative)
subalgebra of $A$ generated by $B$.

\begin{pro}\label{pr:DtinQnF}
There is no unital subalgebra $B$ of the Jordan superalgebra
$Q_n(F)^+$ ($n\geq 2$), isomorphic to $D_t$ ($t\ne 0$), and with
$B'=Q_n(F)$.
\end{pro}
\begin{proof}
Write $A=Q_n(F)$, and take a basis $\{ e,f,u,v\}$ of $B\cong D_t$ as
in Section 1. Since $B$ is a unital subalgebra, $e+f=1_A$.
Therefore, as $e^2=e$, $f^2=f$ and $ef=fe=(1_A-e)e=0$, we may assume
also that
$$e= \begin{pmatrix} I_s & 0 &0 &0 \\
0 & 0 &0 &0\\
0 & 0 &I_s &0\\
0 & 0 &0 &0
\end{pmatrix} ,\quad f= \begin{pmatrix} 0& 0 &0 &0 \\
0 & I_m &0 &0\\
0 & 0 &0&0\\
0 & 0 &0 &I_m
\end{pmatrix}.$$

Consider the Peirce decomposition associated to the idempotents $e$
and  $f$, and note  that $u,v \in A_{\bar{1}} \cap
(Q_n(F)^+)_{1/2}(e) \cap (Q_n(F)^+)_{1/2}(f)$. Hence
$$u= \begin{pmatrix} 0 & 0 &0 &a \\
0 & 0 &b &0\\
0 & a &0 &0\\
b & 0 &0 &0
\end{pmatrix} \mbox{ and } \thickspace v= \begin{pmatrix} 0& 0 &0 &c \\
0 & 0 &d &0\\
0 & c &0&0\\
d & 0 &0 &0
\end{pmatrix},$$
for some $a,c\in M_{s\times m}(F)$, $b,d\in M_{m\times s}(F)$.
 But this contradicts that $B ^\prime$ be equal to $A$, because, for
instance,
$$
\begin{pmatrix} 0 & 0 &x &0 \\
0 & 0 &0 &0\\
x & 0 &0 &0\\
0 & 0 &0 &0
\end{pmatrix} \notin B^{\prime}, \thickspace \mbox{ for } 0\ne x \in M_{s\times s}(F).
$$
This finishes the proof.
\end{proof}

\bigskip

Now, if $Q_n(F)$ is replaced by $M_{p,q}(F)$, some knowledge of the
universal enveloping algebra of $D_t$ is needed.

I. P. Shestakov determined $U(D_t)$ (see \cite {Ma-Ze2}), which is
intimately related to the orthosymplectic Lie superalgebra
$osp(1,2),$ that is, the superalgebra whose elements are the
skew\-symme\-tric matrices of $M_{1,2}(F)$ relative to the
orthosymplectic superinvolution, with Lie bracket  $[a,b]=
ab-(-1)^{\bar a \bar b}ba$:
$$
osp(1,2) =\left\{ \begin{pmatrix} 0 & \beta &\alpha \\ -\alpha & \gamma & \mu
\\ \beta  & \nu & -\gamma \end{pmatrix} \ : \ \alpha, \beta, \mu, \gamma , \nu \in
F \right\}.
$$
The following elements in $osp(1,2)$, which form a basis, will be
considered throughout:
$$h=  \begin{pmatrix}  0 & 0 & 0 \\ 0 & 1 & 0 \\ 0 & 0 & -1
\end{pmatrix}, \thickspace  e= \begin {pmatrix}  0 & 0 & 0 \\ 0 & 0 & 1 \\ 0 & 0 & 0
\end{pmatrix}, \thickspace f= \begin {pmatrix}  0 & 0 & 0 \\ 0 & 0 & 0 \\ 0 & 1 & 0
\end{pmatrix},$$

 $$x=\begin {pmatrix}  0 & 0 & -1 \\ 1 & 0 & 0 \\ 0 & 0 & 0
\end{pmatrix}, \thickspace y= \begin {pmatrix} 0 & 1 & 0 \\ 0 & 0 & 0 \\ 1 & 0 & 0
\end{pmatrix},$$
which verify $[h,e]=2e, \quad [h,f]=-2f, \quad [h,x]=x, \quad
[h,y]=-y, \quad [e,y]=x, \quad [f,x]=y, \quad [x,x]=-2e, \quad
[y,y]=2f, \quad [x,y]=xy+yx=h.$

\bigskip

Then $U(D_t)$ is given by

\begin{theo}\label{th:Ivan}\emph{(I. Shestakov)}
If $t\neq 0, \pm 1,$ then the universal associative enveloping of
$D_t$ is $(U(D_t),\iota)$ where $U(D_t) = U(osp(1,2)) /\ideal
\bigl\langle (xy-yx)^2+(xy-yx) + \frac{t}{(1+t)^2} \bigr\rangle$ and
\begin{alignat*}{3}
\iota: \thickspace &D_t  \quad & \longrightarrow & \quad U(D_t)\\
       &e  \quad &\longmapsto & \quad \iota (e)= \frac{1}{t-1}(t1+ (1+t)\overline{(xy-yx)}),\\
       &f  \quad &\longmapsto & \quad \iota (f)= \frac{1}{1-t}(1+ (1+t)\overline{(xy-yx)}),\\
&u  \quad &\longmapsto & \quad \iota (u)= 2\bar{x},\\
&v  \quad &\longmapsto & \quad \iota (v)= -(1+t)\bar{y},
\end{alignat*}
\noindent  where $\bar {z}$ denotes the class of $z\in osp(1,2)$
modulo the ideal generated by $ (xy-yx)^2+(xy-yx) +
\frac{t}{(1+t)^2} $.
\end{theo}

\bigskip

Here  $U(osp(1,2))$ denotes the universal enveloping algebra of the
Lie superalgebra $osp(1,2)$ (see \cite[section 1.1.3]{KacLie}).

Note that the element $a=\overline{xy-yx}\in U(D_t)$ satisfies
$a^2+a+\frac{t} {(1+t)^2}=0,$ hence if $a^\prime = -(1+t)a$,
${a^\prime}^2-(1+t)a^\prime +t =0$ and in this way the original
version of  Shestakov's Theorem is recovered.

\bigskip

The even part of $osp(1,2)$, which is the span of the elements
$h,e,f$ above, is isomorphic to the three dimensional simple Lie
algebra $ sl(2,F)$, so given  any finite dimensional irreducible
$U(osp(1,2))$-module $V$, by restriction $V$ is also a module for
$sl(2,F)$. The well-known representation theory of $sl(2,F)$ shows
that  $h$ acts diagonally on $V$ (see \cite[7.2 Corollary]{Hum}),
its eigenvalues constitute a sequence of integers, symmetric
relative to $0$, and hence $V$ is the direct sum of the subspaces
$V_m =\{ v\in V : h \cdot v=mv\}$ with $m\in \bZ$.

By finite dimensionality, there exists a largest nonnegative integer
$m$ with  $V_m\neq 0$. Pick a nonzero element
 $v\in V_m$ (a highest weight vector). Changing the parity in $V$ if necessary,
this element $v$ can be assumed to be even.

Since $h(ev)=[h,e]v+e(hv)=(m+2)ev$, it follows that $ev=0$, and
since $h(xv)=[h,x]v+x(hv)=(m+1)xv$, it follows that $xv=0$ too. Let
$\frg=osp(1,2)$, then $\frg=\frg_{-}\oplus\frh\oplus\frg_{+}$, where
$\frg_+=Fe+Fx$, $\frh=Fh$, and $\frg_-=Ff+Fy$, and let $W=W_0=Fw$ be
the module over $\frh+\frg_+$ given by $hw=mw$, $ew=0$, and $xw=0$.
The  map $W \longrightarrow V$ such that $\lambda w \longmapsto
\lambda v$ for any $\lambda \in F$ is a homomorphism of
$(\frh+\frg_+)$-modules, which can be extended to a homomorphism of
$\frg$-modules (that is, of $U(osp(1,2))$-modules) as follows:
$$
\begin{array}{ccc} \varphi : U(\frg)\otimes _{U(\frh+\frg_+)} W
    & \longrightarrow & V \\
a\otimes w & \longmapsto & av. \end{array}$$

\noindent Since $V$ is and irreducible $osp(1,2)$-module, $\varphi$
is onto. We denote by $U(m)$ the $U(\frg)$-module $U(\frg)\otimes
_{U(\frh+\frg_+)} W$ and identify the element $1\otimes w$ with $w$.
Then:

$$\begin{array}{rclrcl} hy^iw & = & (m-i)y^iw, & fy^iw & = & y^{i+2}w, \\
xy^{2i}w & = & -iy^{2i-1}w,  & xy^{2i+1}w & = & (m-i)y^{2i}w, \\
ey^{2i}w & = & i(m - i+1)y^{2i-2}w, & ey^{2i+1}w & = &
i(m-i)y^{2i-1}w,
\end{array}$$

\noindent and hence it follows that the set $\{ w, yw, y^2w, \dots
\}$ spans the vector space $U(m)$. We remark that $I_m=\mathop{span}
\langle y^{2m+1}w, y^{2m+2}w, \dots \rangle$ is a proper submodule
of $U(m)$, and because $V$ is irreducible and the weights of the
elements $y^{2m+i}w$ are all different from $m$, it follows that
$\varphi(I_m)\ne V$, so by irreducibility $\varphi(I_m)=0$. Thus the
set $\{v, yv, y^2v, \dots , y^{2m}v\}$ spans the vector space $V$.
Again, the theory of modules for $sl(2,F)$ shows that
$v,y^2v,\ldots,y^{2m}v$ are all nonzero (see \cite[7.2]{Hum}), and
hence so are the elements   $yv, y^3v, \dots, y^{2m-1}v$. Note that
the elements $v,yv, y^2v, \dots , y^{2m}v$ are linearly independent,
as they belong to different eigenspaces relative to the action of
$h$. We conclude that $\{ v,yv, y^2v, \dots, y^{2m}v\}$ is a basis
of $V$.

Denote $V$ by $V(m)$ and write $e_i=y^iv$. Then,
\begin{alignat*}{2}
V(m)_{\bar{0}} \quad &=  \quad \langle e_0, e_2, \ldots , e_{2m}
 \rangle \thickspace ,\\
V(m)_{\bar{1}} \quad  &=  \quad \langle e_1, e_3, \ldots , e_{2m-1}
\rangle \thickspace .
\end{alignat*}

Observe that
\begin{alignat*}{3}
(xy-yx)e_{2i} \hspace{0.4cm} \quad &= \quad \hspace{0.5cm}
  (m-i)e_{2i}+ ie_{2i} \quad &=&  \qquad \hspace{0.25cm} me_{2i} \thickspace ,\\
(xy-yx)e_{2i+1} \quad &= \quad xe_{2i+2}- (m-i)e_{2i+1} \quad &=&
\quad -(m+1)e_{2i+1} \thickspace ,
\end{alignat*}

\noindent and so the minimal polynomial of the action of $xy-yx$ is
$(X-m)(X+(m+1)) = X^2+X-m(m+1)$, and therefore the finite
dimensional irreducible $U(osp(1,2))$-modules coincide with the
irreducible modules for $U(osp(1,2)) / \ideal \langle (xy-yx)^2
+(xy-yx) -m(m+1) \rangle$.

\bigskip

Therefore, if $V$ is a finite dimensional irreducible
$U(D_t)$-module ($t\ne 0,\pm 1$), then by Shestakov's Theorem
(Theorem \ref{th:Ivan}), $V$ is an irreducible module for $osp(1,2)$
in which the minimal polynomial of the action of $xy-yx$ divides
$X^2+X+\frac{t}{(1+t)^2}$. From our above discussion, there must
exist a natural number $m$ such that $\frac{t}{(1+t)^2}=-m(m+1)$,
that is, either $t=-\frac{m}{m+1}$ or $t=-\frac{m+1}{m}$. Thus,

\bigskip

\begin{coro} (C. Mart\'{\i}nez, E. Zelmanov)\quad
The universal enveloping algebra $U(D_t)$ ($t\ne 0,\pm 1$) has a
finite dimensional irreducible module if and only if there exists a
natural number $m$ such that either $t= - \frac{m}{m+1}$ or $t= -
\frac{m+1}{m}.$ In this case, up to parity exchange, its unique
irreducible module is $V(m)$ (that is, the irreducible module for
$U(osp(1,2))$ annihilated by the ideal generated by $
(xy-yx)^2+(xy-yx)-m(m+1)$).
\end{coro}

\bigskip

Something can be added here:

\bigskip

\begin{pro}\label{pr:superform} Up to scalars, the module
$V(m)$ has a unique nonzero even bilinear
form $(.\mid .)$ such that $\rho_x$ and $\rho_y$, the multiplication
operators by $x$ and $y$, are supersymmetric (that is,
$(zv|w)=(-1)^{|v|} (v|zw)$ for any $v,w \in V_{\bar{0}}\cup
V_{\bar{1}}$ with $z=x,y$).
\end{pro}
\begin{proof}
 If $\rho_x, \rho_y$ are supersymmetric then
$\rho_{[x,x]}=2\rho_x^2$, $\rho_{[y,y]}=2\rho_y^2$, and
$\rho_{[x,y]}=\rho_x\rho_y+\rho_y\rho_x$ are skewsymmetric, that is,
$\rho_e$, $\rho_f$, and $\rho_h$ are skewsymmetric. But $\rho_h$
being skewsymmetric implies that $(V_{(\alpha)}| V_{(\beta)})=0$ if
$\alpha +\beta \neq 0$, where  $V_{(\alpha)} =\{ v\in V(m) \ : \
hv=\alpha v \}$, because $(hV_{(\alpha)}| V_{(\beta)})=
-(V_{(\alpha)}|hV_{(\beta)})$, and therefore $(\alpha
+\beta)(V_{(\alpha)}| V_{(\beta)})=0$.
 Hence we can check that $(.\mid
. )$ is determined by $(e_0|e_{2m})$, as

$$(e_1|e_{2m-1})=(ye_0|e_{2m-1})= (e_0|ye_{2m-1})=(e_0|e_{2m}).$$

So, up to scalars, it can be assumed that  $(e_0|e_{2m})=1$.

Using  that $\rho_y$ is supersymmetric , recursively we get

$$
 \begin {array}{l} (e_{2r}|e_{2(m-r)})=(-1)^r, \\ (e_{2r+1}|e_{2(m-r)-1})=(-1)^r
\end{array}
$$

 \noindent and  $(e_i|e_j)=0$ otherwise. Now it can be checked that $\rho_x$ is
 supersymmetric too.
 \end{proof}

\bigskip

Note that $(.\mid . )$ is supersymmetric if $m$ is even and
superskewsymmetric if $m$ is odd. In the latter case, one can
consider $V(m)^{op}$ with the supersymmetric bilinear superform
given by $(u|v)^\prime = (-1)^{|u|}(u|v)$ where $|u|$ denotes the
parity in $V(m)$.

\bigskip

Consider again the finite dimensional irreducible $U(D_t)$-module
($t=-\frac{m}{m+1}$ or $t=-\frac{m+1}{m}$) $V=V(m)$, with the
bilinear superform in the proposition above. It is known that this
determines a superinvolution in $A=\End_F(V)$ such that every
homogeneous element $f\in \End_F(V)$ is mapped to $f^*$
 verifying $(fv,w)=(-1)^{\bar f\bar v}(v,f^*w)$. Note that, since
 $\rho_x$ and $\rho_y$ are supersymmetric,
$D_t$ is thus embedded in  $H(\End_F(V),*)$ as follows:

\begin{eqnarray*}
D_t &\longrightarrow & H(\End_F(V),*) \\
 e & \longmapsto & \frac{1}{t-1} (t\rho_{Id}+ (1+t)(\rho_x \rho_y - \rho_y \rho_x))\\
 f & \longmapsto & \frac{1}{1-t} (\rho_{Id}+ (1+t)(\rho_x \rho_y - \rho_y \rho_x))\\
 u  & \longmapsto & 2 \rho_x \\
 v &\longmapsto & -(1+t)\rho_y.
\end{eqnarray*}
 Moreover, unless $t \neq -2, -1/2$ (that is, unless $m=1$), by dimension count,
one has $D_t \subsetneqq H(\End_F(V),*)$.

\bigskip

The conclusion of all these arguments is the following:

\begin{pro}\label{pr:DtEndV}
Let $V$ be a nontrivial finite dimensional vector superspace and let
$B$ be a unital subalgebra of the simple Jordan superalgebra
$\End_F(V)^+$, isomorphic to $D_t$ ($t\ne 0,\pm 1$), and such that
$B'=\End_F(V)$. Then one of the following situations holds:

\begin{enumerate}

\item[{\rm (i)}] either $t= -\frac{m}{m+1}$ or $t= -\frac{m+1}{m}$ for an even number
 $m$,  such that $V\cong V(m)$, and through this isomorphism $B
\subseteq H(\End_F(V), *)$ where $*$ is the superinvolution
associated to the bilinear superform of Proposition
\ref{pr:superform},

\item[{\rm (ii)}] or $t=-\frac{m}{m+1}$ or $t=-\frac{m+1}{m}$
for an odd number $m$  such that $V\cong V(m)^{op}$ and through this
isomorphism $D_t \subseteq H(\End_F(V), \diamond)$, where $\diamond$
is the superinvolution associated to the bilinear superform $(.\mid
.)^\prime.$

\end{enumerate}

\end{pro}
\begin{proof} The hypotheses imply that there is a surjective
homomorphism of associative algebra $U(D_t)\rightarrow \End_F(V)$,
so $V$ becomes an irreducible module for $U(D_t)$ and the arguments
above apply.
\end{proof}

\bigskip

Since the superalgebra $\End_F(V)$, for a superspace $V$, is
isomorphic to $M_{p,q}(F)$, for $p=\dim V\subo$, $q=\dim V\subuno$,
the next result follows:

\begin{coro}\label{co:DtinMpqF}
The simple Jordan superalgebra $M_{p,q}(F)^+$ contains a unital
 subalgebra $B$, isomorphic to $D_t$ ($t\ne 0,\pm 1$), and such that $B'=M_{p,q}(F)$,
 if and only if $q=p\pm 1$ and either $t=-\frac{p}{q}$, or $t=-\frac{q}{p}$.
\end{coro}

\bigskip

Proposition \ref{pr:DtinQnF} and Corollary \ref{co:DtinMpqF} give
all the possibilities for embeddings of the Jordan superalgebra
$D_t$ ($t\ne 0,\pm 1$) as unital subalgebras in $A^+$, in such a way
that the associative subalgebra generated by $D_t$ is the whole $A$,
$A$ being a simple associative superalgebra. For these cases, one
always has
 $D_t\subseteq H(A,*)$, for a suitable superinvolution. By dimension
count, equality is only possible here if $t=-2$ (or
$t=-\frac{1}{2}$). This corresponds to the isomorphism $D_{-2}\cong
osp_{1,2}$.

\bigskip

For later use, let us recall the following results on universal
enveloping algebras of some other Jordan superalgebras (see
\cite{Ma-Ze2}):

\begin{theo}\label{th:p2} \emph{(C.~Mart\'{\i}nez and E.~Zelmanov)}
\begin{enumerate}
\item[{\rm (i)}] The universal enveloping algebra of $ p(2) $ is isomorphic to  $M_{2,2}(F[t]),$
 where $F[t]$ is the polynomial algebra in the variable $t$.
\item[{\rm (ii)}] The universal enveloping algebra of $M_{1,1}(F)$ is $(U(D), u)$ with
\end{enumerate}
\begin{multline*}
U(D)= \begin{pmatrix} F[z_1, z_2]+ F[z_1, z_2]a & 0 \\
0 & F[z_1, z_2]+ F[z_1, z_2]a
\end{pmatrix} \\
\oplus \begin{pmatrix} 0& F[z_1, z_2]+ F[z_1, z_2]a^{-1}z_2  \\
F[z_1, z_2]z_1+ F[z_1, z_2]a & 0
\end{pmatrix}
\end{multline*}
 where $z_1, z_2$ are variables,  $a$ is a root of $X^2 +X -z_1 z_2 \in
 F[z_1, z_2]$, and $u: M_{1,1}(F) \rightarrow U(D)^+$ is given by
$$\begin{pmatrix} \alpha_{11} & \alpha_{12} \\
\alpha_{21} & \alpha_{22}
\end{pmatrix} \mapsto \begin{pmatrix} \alpha_{11} & \alpha_{12}+ \alpha_{21}a^{-1}z_2  \\
\alpha_{12}z_1+ \alpha_{21}a & \alpha_{22}
\end{pmatrix}.$$
\end{theo}

\bigskip

\begin{theo}\label{th:others} \emph{(C.~Mart\'{\i}nez and E.~Zelmanov)}
\begin{enumerate}
\item[{\rm (i)}] $U(M_{m,n}^{+})(F) \cong M_{m,n}(F) \oplus M_{m,n}(F)$ for $(m,n) \neq (1,1);$
\item[{\rm (ii)}] $U(Q^{+}_n(F))= Q_n(F) \oplus Q_n(F),$ $n \geq 2 ;$
\item[{\rm (iii)}] $U(osp_{m,n}(F)) \cong M_{m,n}(F),$ $(m,n) \neq (1,2);$
\item[{\rm (iv)}]$U(p(n)) \cong M_{n,n}(F),$ $n \geq 3.$
\end{enumerate}

\end{theo}

\bigskip

\section {Maximal subalgebras of $A^+.$}

Let $B$ be a maximal subalgebra of $A^+$, $A$ being a simple
associative superalgebra (so $A$ is isomorphic to either
$M_{p,q}(F)$ or $Q_n(F)$, for some $p$ and $q$, or $n$). If
$B^\prime \neq A$ then $B^\prime \subseteq C$ with $C$ a maximal
subalgebra of the associative superalgebra $A$, and then $C^+ =B$ by
maximality. Therefore a maximal subalgebra of $A^+$ is of one of the
following types, either:

\begin{enumerate}

\item[{\rm (i)}] $B^{\prime}=A$ and $B$ is semisimple, or

\item[{\rm (ii)}]  $B=C^+$ with $C$ a maximal subalgebra of $A$ as associative
superalgebra, or

\item[{\rm (iii)}]  $B^{\prime}=A$ and $B$ is not semisimple.

\end{enumerate}

\bigskip

\subsection{$B^\prime = A$ and $B$ semisimple.}

Let us assume first that $B$ is a maximal subalgebra of the simple
superalgebra $A^+$, with $B'=A$ and $B$ semisimple.

For the moment being, let us drop the maximality condition, so let
us suppose that $B$ is just a semisimple subalgebra of $A^+$ with
$B'=A$. By Theorem \ref{th:semisimple}, $B=\sum_{i=1}^{r}(
J_{i1}\oplus \dots \oplus J_{ir_i} + Fe_i)\oplus M_1\oplus \dots
\oplus M_t$ where $M_1, \dots, M_t$ are simple Jordan superalgebras
and $J_{ij}$ are Kaplansky superalgebras.

We claim that $B$ has neither direct summands $M_i$ isomorphic to
the Kaplansky superalgebra $K_3$ nor direct summands of the type
$(J_{i1}\oplus \dots \oplus J_{ir_i} + Fe_i)$. Indeed, otherwise
$A^+$ would contain a subalgebra isomorphic to $K_3$. Let $e$ be its
nonzero even idempotent and $x,y$ odd elements with $x\cdot y=e$.
Then, in the associative superalgebra $A$ (which is isomorphic to
either $M_{p,q}(F)$ or $Q_n(F)$, and hence there is a trace form),
one has $\trace(e)=\trace(x\cdot y)=\frac{1}{2}\trace(xy-yx)=0$.
However, any nonzero idempotent in a matrix algebra over a field of
characteristic $0$ has nontrivial trace. A contradiction.

Therefore,
 $B = M_1 \oplus \dots \oplus M_t$, where the $M_i$'s are unital simple Jordan
superalgebras.

Consider now the identity element $f_i$  of each $M_i.$ Then
$B=f_1Bf_1 \oplus \ldots \oplus f_tBf_t$. If $t>1$, it follows that
$B^\prime \subset f_1Af_1 \oplus (1-f_1)A(1-f_1)\subsetneqq A$, a
contradiction. Hence  $B$ is simple and, therefore,  is isomorphic
to one of the following special superalgebras: $ D_t$, $H(D,*)$ (for
a simple associative superalgebra $D$ with superinvolution $*$), the
superalgebra of a superform, or $D^+$ for a simple associative
superalgebra $D$. (Recall that $K_{10}$ and the Kantor superalgebra
are exceptional superalgebras.)

In case  $B$ were the superalgebra of a superform over a vector
superspace $V$, let $x,y\in V_{\bar{1}}$ such that $x \cdot y =1_A$.
Then $x \cdot y = \frac{1}{2}(xy-yx)=1_A$, and  again $\trace (x
\cdot y) =0\neq \trace (1_A)$, a contradiction that shows that
$V_{\bar{1}}=0$. But then $B\subseteq A_{\bar{0}}$ and $B^\prime
\subseteq A_{\bar{0}} \neq A$, contrary to our hypotheses.

Now,  in case $B$ is isomorphic to $D_t$ ($t\ne 0$), Proposition
\ref{pr:DtinQnF} shows that $A$ is not isomorphic to $Q_n(F)$ and
Corollary \ref{co:DtinMpqF} shows that $B$ is never a maximal
subalgebra of $A\cong M_{p,q}(F)$ unless $t=-2$ (or $-\frac{1}{2}$).
In this case $B$ is isomorphic to $H(D,*)$ for a suitable $(D,*)$.

Therefore:

\begin{lemma}\label{le:BinAplus}
Let $B$ be a subalgebra of the Jordan superalgebra $A^+$, where $A$
is a finite dimensional simple associative superalgebra over an
algebraically closed field $F$ of characteristic $0$. If
$B^{\prime}=A$ and $B$ is semisimple, then either $B$ is isomorphic
to $D_t$ ($t\ne 0,1,-1,-2,-\frac{1}{2}$), or $B= D^+$ or $B=H(D,*)$,
for  a simple associative superalgebra $D$ and a superinvolution
$*$. Moreover, if $B$ is a maximal subalgebra of $A^+$, then  the
first possibility does not hold.
\end{lemma}

\bigskip

Our next goal consists in proving that, in case $B=D^+$ or $B=
H(D,*)$, one has that $D$ is isomorphic to $A$. For this  the
following result (see \cite {GoT}) will be used:

\begin{theo}\label{th:Carlos} \emph{(C. G\'omez-Ambrosi)}
Let $S$ be a unital associative superalgebra with superinvolution
$*$. Assume that the following conditions hold:

\begin{enumerate}

\item[{\rm (i)}] $S$ has at least three symmetric orthogonal idempotents.

\item[{\rm (ii)}] If $S=\sum_{i=1}^n S_{ij}$ is the Peirce decomposition related to
them, then $S_{ij}S_{ji}=S_{ii}$ holds for $i,j=1, \dots ,n$,
\end{enumerate}

\noindent and let $\phi \colon H(S,*) \to (A, \cdot)^+$ be a
homomorphism of Jordan superalgebras, for an associative
superalgebra $(A,\cdot)$. Then $\phi$ can be extended univocally to
an associative homomorphism $\varphi \colon S \to A$.

\end {theo}

\bigskip

We shall proceed in several steps, where the assumptions are that
$B$ is just a semisimple subalgebra of $A^+$ with $B'=A$:

\medskip

\noindent\textbf{a)}\quad Assume first that $B=H(D,*)$ for a simple
associative superalgebra with involution $(D,*)$. Let us denote the
multiplication in $D$ by $\diamond$. The inclusion map $\iota\colon
B=H(D,*) \to (A,\cdot)^+$ is a Jordan homomorphism. Then (Section
1), $D$ is isomorphic to $M_{p,q}(F)$, for suitable $p,q$, and $*$
corresponds to either the transpose involution or an orthosymplectic
involution. If neither $D$ is a quaternion superalgebra (isomorphic
to $M_{1,1}(F)$), nor $H(D,*)$ is isomorphic to $p(2)$ or
$osp_{1,2}$, then $D$ satisfies the hypotheses of Theorem
\ref{th:Carlos} and, therefore, $\iota:B\rightarrow A$ can be
extended to an associative homomorphism $\tau:D\rightarrow A$. But
the subalgebra $B'$ generated by $B$ in $A$ is the whole $A$. Hence
$\tau$ is onto and, as $D$ is simple, it is one-to-one too.
Therefore $D$ is isomorphic to $A$. Thus, we are left with three
cases:

\medskip

\noindent\textbf{a.1)}\quad If $H(D,*)$ is isomorphic to $osp_{1,2}$
then, since $osp_{1,2}$ is isomorphic to $D_{-2}$, $H(D,*)$ is
isomorphic to $D_{-2}$.

\medskip

\noindent\textbf{a.2)}\quad If $D$, with multiplication $\diamond$,
is isomorphic to $M_{1,1}(F)^+$, with superinvolution $*$ as in 6)i)
in Section 1, then $H(D,*)$ is isomorphic to $F1+Fu$, with $u^2=0$.
Thus, the universal enveloping algebra of $H(D,*)$ is $F[u]$, the
ring of polynomials over $F$ on the variable $u$, and there exists
an associative homomorphism $\varphi: F[u]\rightarrow A$, which
extends $\iota:B\rightarrow A$. Again, $\varphi$ is onto since
$B'=A$. Therefore $A$ should be commutative, a contradiction.

\medskip

\noindent\textbf{a.3)}\quad Finally, if $H(D,*)$ is isomorphic to
$p(2)$, Theorem \ref{th:p2} shows that its universal enveloping
algebra is isomorphic to $M_{2,2}(F[t])$, where $F[t]$ is the
polynomial algebra on the indeterminate $t$. As before, this gives a
surjective homomorphism $\phi:M_{2,2}(F[t])\rightarrow A$. Recall
that $A$ is isomorphic either to $M_{p,q}(F)$ or to
$Q_n(F)=M_n(F)\oplus M_n(F)u$ ($u^2=1$). Let $e_1,e_2,e_3,e_4$ be
primitive orthogonal idempotents of $M_{2,2}(F)$, with $e_1+e_2$ and
$e_3+e_4$ being the unital elements in the two simple direct
summands of the even part. Since the restriction of $\phi$ to
$M_{2,2}(F)$ is injective because $M_{2,2}(F)$ is simple, the images
$\phi(e_1),\phi(e_2),\phi(e_3),\phi(e_4)$ are nonzero orthogonal
idempotents in $A\subo$ with $\sum_{i=1}^4\phi(e_i)=1_A$. Write
$U=M_{2,2}(F[t])$ and
 consider
the Peirce decomposition of $U$ relative  to $e_1,e_2,e_3,e_4,$: $U=
\sum U_{ij}$, and the Peirce decomposition of $A$ relative  to $\phi
(e_1),\phi(e_2), \phi(e_3), \phi(e_4)$: $A=\sum A_{ij}$. Since
$U_{ii}$ is isomorphic to  $F[t]$, it follows that $A_{ii}$ is
commutative (as a quotient of $F[t]$) for any $i=1,2,3,4$. Therefore
either $p+q=4$  or $n=4$, that is $A\cong Q_4(F)$. Consider now the
restriction $\phi|_{M_{2,2}(F[t])\subo}
 \colon M_{2,2}(F[t])\subo \to
A$. If $A\cong M_{p,q}(F),$ with $p+q=4$ one has that
$\phi(M_{2,2}(F[t])\subo)
 = \phi(M_2(F[t])) \oplus
\phi(M_2(F[t])) = A_{\bar{0}} \cong M_p(F) \oplus M_q(F)$, and
therefore $p=2$ and $q=2$, and $D \cong M_{2,2}(F) = A$. If $A\cong
Q_4(F)$, then $(M_2(F[t])\times {0})$ is an ideal of
$M_{2,2}(F[t])_{\bar{0}}$, and so $\phi (M_2(F[t])\times {0})$ is an
ideal of $A_{\bar{0}} \cong M_4(F)$. Since $M_4(F)$ is simple and
$\phi (e_1), \phi(e_2)$ are nonzero idempotents, it follows that
$\phi (M_2(F[t])\times {0})=A_{\bar{0}},$ and so $\phi (e_1)+
\phi(e_2) =1_A,$ that is a contradiction because $\phi
(e_1)+\phi(e_2)+\phi(e_3)+\phi(e_4) =1$, with $\phi(e_3), \phi(e_4)$
nonzero orthogonal idempotents.

\bigskip

\noindent\textbf{b)}\quad Assume now that $B=D^+$ for a simple
associative superalgebra $D$. Consider the opposite superalgebra
$D^{op}$ defined on the same vector space as $D$, but with the
multiplication given by $a\diamond b = (-1)^{\bar a \bar b} b \cdot
a$, and the direct sum $D\oplus D^{op}$, which is endowed with the
superinvolution $- \colon D\oplus D^{op} \to D\oplus D^{op}$, such
that $\overline {(x,a)} = (a,x)$. Note that if $e_1, e_2,\dots, e_n$
are orthogonal idempotents in $D$, then $(e_1,e_1), (e_2,e_2),
\dots, (e_n,e_n)$ are also orthogonal idempotents in $D\oplus
D^{op}$, and the Peirce spaces are given by $(D\oplus D^{op})_{ij} =
D_{ij}\oplus (D^{op})_{ji}$. So if $D$ satisfies conditions (i) and
(ii) in  Theorem \ref{th:Carlos}, then so does $D\oplus D^{op}$.
Since $D^+$ is isomorphic to $H(D\oplus D^{op}, -)$, there is a
homomorphism of Jordan superalgebras $\phi \colon H(D\oplus
D^{op},-) \to A^+$.

\medskip

\noindent\textbf{b.1)}\quad Suppose that $D$ is not isomorphic
 to $M_{1,1}(F),$ nor to $Q_2(F)$, then from Theorem
\ref{th:Carlos}, $\phi$ can be extended to an associative
homomorphism $\varphi \colon D\oplus D^{op} \to A$.
 As before,
$\varphi$ is onto because $B^\prime =A,$ so $D\oplus D^{op} /Ker
\varphi$ is isomorphic to $A$ and either $Ker \varphi \cong D$ or
$Ker \varphi \cong D^{op}$, because $A$ is simple. Hence either
$D\cong A$ or $D^{op}\cong A$, that is, $\dim D= \dim A$, a
contradiction.

\medskip

\noindent\textbf{b.2)}\quad If $D$ is isomorphic to $M_{1,1}(F)$
(that is, $D$ is a quaternion superalgebra), consider the universal
enveloping algebra $(U(D),u)$ of $D^+$ (see Theorem \ref{th:p2}).
The Jordan homomorphism $\iota\colon D \to A^+$ extends to  an
associative homomorphism $\varphi \colon U(D) \to A$ such that
$\varphi \circ u =\iota$. But $B^\prime =A,$ and hence it follows
that $\varphi$ is onto and, therefore, $U(D)/Ker \varphi \cong A$.
Recall that $F$, the underground field, is assumed to be
algebraically closed, so either $A\cong Q_n(F)$ or $A\cong
M_{p,q}(F).$ But $(U(D)/Ker \varphi)_{\bar{0}}$ is commutative, so
$A_{\bar{0}}$ is commutative and therefore either $A\cong Q_1(F)$ or
$A\cong M_{1,1}(F)$, a contradiction to $D$ being isomorphic to
$M_{1,1}(F)$.

\medskip

\noindent\textbf{b.3)}\quad Otherwise $D$ is isomorphic to $Q_2(F)$,
and hence the universal enveloping algebra $(U(D),u)$ of $D^+$ is
isomorphic to $D\oplus D$ (see Theorem \ref{th:others}). Hence there
is a surjective homomorphism $\varphi \colon U(D) \to A$ which
extends $\iota$. As before, $\varphi$ is onto and so $U(D)/ Ker
\varphi \cong A.$ But $A$ is simple, so $Ker \varphi \cong D $ and
$A\cong D$, a contradiction.

\bigskip

Therefore, Lemma \ref{le:BinAplus} can be improved to:

\begin{lemma}\label{le:BinAplusimproved}
Let $A$ be a finite dimensional simple associative superalgebra over
$F$, and let $B$ be a semisimple subalgebra of $A^+$  with $B'=A$,
then either $B$ is isomorphic to $D_t$ ($t\ne 0,\pm
1,-2,-\frac{1}{2}$), or $B$ equals $H(A,*)$, for a superinvolution
$*$. Moreover, if $B$ is a maximal subalgebra of $A^+$, then
$B=H(A,*)$ for a superinvolution $*$ of $A$.
\end{lemma}

In consequence, if $B$ is a maximal subalgebra of $A$, which is
semisimple and satisfies $B^\prime =A$, Lemma
\ref{le:BinAplusimproved} shows that $B$ coincides with the
subalgebra of hermitian elements of $A$ relative to a suitable
superinvolution. The converse also holds:

\begin{theo}\label{th:BsemiAplus}
Let $A$ be a finite dimensional simple associative superalgebra over
an algebraically closed field of characteristic zero, and let $B$ be
a semisimple subalgebra of $A^+$ such that $B^\prime =A $. Then $B$
is a maximal subalgebra of $A^+$ if and only if there is a
superinvolution $*$ in $A$ such that  $B= H(A,*)$.
\end{theo}
\begin{proof}
The only thing left is to show that if $A$ is a finite dimensional
simple associative superalgebra endowed with a superinvolution $*$,
then $H(A,*)$ is a maximal subalgebra of $A^+$.

Our hypotheses on the ground field imply that, up to isomorphism, we
are left with the next two  possibilities:

\begin{enumerate}
\item[i)] $A=M_{n,n}(F)$, and $\thickspace \begin{pmatrix} a & b \\
c & d \end{pmatrix}^{*} = \begin{pmatrix} d^t & -b^t \\
c^t & a^t \end{pmatrix}.$
\item[ii)] $A=M_{n,2m}(F)$, and $\thickspace \begin{pmatrix} a & b \\
c & d \end{pmatrix}^* = \begin{pmatrix} a^t & c^tq \\
-q^tb^t & q^td^tq \end{pmatrix},$ where $q= \begin{pmatrix} 0 & I_m \\
-I_m & 0 \end{pmatrix}.$
\end{enumerate}

Note that $A=H\oplus K$, where $H=H(A,*)$ and $K$ is the set of
skewsymmetric elements of $(A,*)$.

\noindent\textbf{i)}\quad In the first case
\[
H= \left\{ \begin{pmatrix} a & b \\
c & a^t \end{pmatrix} : c \mbox{ symmetric, } b \mbox{
skewsymmetric } \right\} ,\]
\[ K= \left\{ \begin{pmatrix} a & b \\
c & -a^t \end{pmatrix} : b \mbox{ symmetric, } c \mbox{
skewsymmetric } \right\} ,
\]
and to check that $H(A,*)$ is a maximal subalgebra of $A^+$ it
suffices to prove that $\Jalg\langle H,x\rangle=A^+$ for any nonzero
homogeneous element $x\in K $. ($\Jalg\langle S\rangle$ denotes the
subalgebra generated by $S$.)

If $0\ne x\in K_{\bar{0}}$ then
 $$ x= \begin{pmatrix} a & 0 \\
0 & -a^t \end{pmatrix}  $$

\noindent with $a\in M_n(F)$ and so
$$  \begin{pmatrix} a & 0 \\
0 & -a^t \end{pmatrix} + \begin{pmatrix} a & 0 \\
0 & a^t \end{pmatrix}= \begin{pmatrix} 2a & 0 \\
0 & 0 \end{pmatrix} \in \Jalg \langle H,x \rangle.$$

\noindent We claim that if $\bigl(\begin{smallmatrix}
 a & 0 \\
0 & 0 \end{smallmatrix}\bigr) \in \Jalg \langle H,x \rangle$, then
$\bigl( \begin{smallmatrix} u & 0 \\
0 & 0 \end{smallmatrix}\bigr) \in \Jalg \langle H,x \rangle$, for
any $u\in M_n(F)$. Similarly, if $\bigl(\begin{smallmatrix}
 0 & 0 \\
0 & a \end{smallmatrix}\bigr) \in \Jalg \langle H,x \rangle$, then
$\bigl( \begin{smallmatrix} 0 & 0 \\
0 & u \end{smallmatrix}\bigr) \in \Jalg \langle H,x \rangle$, for
any $u\in M_n(F)$. Actually, since $M_n(F)^+$ is simple and the
ideal generated by $a$ in $M_n(F)^+$ is the vector subspace spanned
by $\{\langle L_{b_1} \ldots  L_{b_m}(a):m\in\bN,\,
b_1,\ldots,b_m\in M_n(F)\}$ ($L_b$ denotes the left multiplication
by $b$ in $M_n(F)^+$), it is enough to realize that
\[
\begin{pmatrix} L_{b_1} \ldots L_{b_m}(a) & 0 \\
0 & 0 \end{pmatrix}=
L_{\bigl( \begin{smallmatrix} b_1 & 0 \\
0 & b_1^t \end{smallmatrix}\bigr) }  \ldots  L_{\bigl( \begin{smallmatrix} b_m & 0 \\
0 & b_m^t \end{smallmatrix}\bigr)} \begin{pmatrix} a  & 0 \\
0 & 0 \end{pmatrix} \in \Jalg \langle H,x \rangle.
\]

So, if $0\ne x \in K_{\bar{0}}$, then $A_{\bar{0}} \subseteq \Jalg
\langle H,x \rangle$. In order to prove that
 $A_{\bar{1}}\subseteq \Jalg \langle H,x \rangle$, note that
$$
\begin{pmatrix} 0 & 0 \\
I_n & 0 \end{pmatrix} \in H,
$$
\noindent and since
$$
\begin{pmatrix} 0 & 0 \\
I_n & 0 \end{pmatrix} \circ
\begin{pmatrix} d & 0 \\
0 & 0 \end{pmatrix} = \frac{1}{2}
\begin{pmatrix} 0 & 0 \\
d & 0 \end{pmatrix}
$$

\noindent it follows that
$$
\begin{pmatrix} 0 & 0 \\
u & 0 \end{pmatrix} \in \Jalg \langle H,x \rangle\quad\text{for any
$u\in M_n(F)$.}
$$

\noindent It remains to prove that
$$
\begin{pmatrix} 0 & u \\
0 & 0 \end{pmatrix} \in \Jalg \langle H,x \rangle\quad\text{for any
$u\in M_n(F)$,}
$$
and the above implies that
$$
\begin{pmatrix} 0 & b \\
0 & 0 \end{pmatrix} \in \Jalg \langle H,x \rangle
$$
for any nonzero skewsymmetric matrix $b$. But
\[
\left( \begin{pmatrix} 0 & b \\
0 & 0 \end{pmatrix} \circ \begin{pmatrix} 0 & 0 \\
0 & M_n(F) \end{pmatrix} \right) \circ \begin{pmatrix} M_n(F) & 0 \\
0 & 0 \end{pmatrix} =
 \begin{pmatrix} 0 & M_n(F)bM_n(F) \\
0 & 0 \end{pmatrix} \subseteq \Jalg \langle H,x \rangle
\]
and $M_n(F)bM_n(F)$ is a nonzero ideal of the simple algebra
$M_n(F)$, so 
it is the whole $M_n(F)$ and
\[
\begin{pmatrix} 0 & M_n(F) \\
0 & 0 \end{pmatrix} \subseteq \Jalg \langle H,x \rangle .
\]

Therefore, $\Jalg \langle H,x \rangle = A^+$ for any nonzero element
$x\in K_{\bar{0}}$.

Now, if $0\ne x\in K_{\bar{1}}$, then
$$ x=\begin{pmatrix} 0 & b \\
c & 0 \end{pmatrix}  $$

\noindent with $b$ a symmetric and $c$  a skewsymmetric $n\times
n$-matrix respectively. Let $y\in H_{\bar{1}},$
$$ y=\begin{pmatrix} 0 & \bar{b} \\
\bar{c} & 0 \end{pmatrix}  $$

\noindent with $\bar{b}$ skewsymmetric and $\bar{c}$ symmetric, such
that $x\circ y \neq 0$. Since $0\ne x\circ y\in K_{\bar{0}}$  we are
back to the `even' case, and so $\Jalg \langle H,x \rangle =A^+$.

\medskip
\noindent\textbf{ii)}\quad In the second case (orthosymplectic
superinvolution), $A=M_{n,2m}(F)$ and
\[
H(A, *)= \left\{ \begin{pmatrix} a & b \\
-q^tb^t & d \end{pmatrix} : a \mbox{ symetric, }   d= \begin{pmatrix} d_{11} & d_{12} \\
d_{21} & d_{11}^t \end{pmatrix}, d_{12}, d_{21} \mbox{ skewsymmetric } \right\},
\]

\[
  K(A, *)= \left\{ \begin{pmatrix} a & b \\
q^tb^t & d \end{pmatrix} : a \mbox{ skewsymmetric, }  d= \begin{pmatrix} d_{11} & d_{12} \\
d_{21} & -d_{11}^t \end{pmatrix}, d_{12}, d_{21} \mbox{ symmetric } \right\}.
\]

\noindent We claim that $\Jalg \langle H,x \rangle =A^+$ for any
nonzero homogeneous element $x\in K$. If $0\ne x\in K_{\bar{1}}$,
then
$$ x=\begin{pmatrix} 0 & b \\
q^tb^t & 0 \end{pmatrix}  $$

\noindent and so
$$ x + \begin{pmatrix} 0 & b \\
-q^tb^t & 0 \end{pmatrix} = \begin{pmatrix} 0 & 2b \\
0 & 0 \end{pmatrix} \in \Jalg \langle H,x \rangle$$

\noindent with $b\in M_{n\times 2m}(F)$. Suppose that
$\bigl(\begin{smallmatrix} 0& b\\ 0&0\end{smallmatrix}\bigr) =
\sum_{i=1, j= n+1}^{n,n+2m} \lambda _{ij}e_{ij}$ with $\lambda =
\lambda_{pq}\neq 0$, where, as usual, $e_{ij}$ denotes the matrix
whose $(i,j)$-entry is $1$ and all the other entries are $0$, then
\[   \left( e_{pp} \circ \begin{pmatrix} 0 & b \\
0 & 0 \end{pmatrix} \right) \circ (e_{qq} + e_{q\pm m, q\pm m})=
\frac{1}{4} (\lambda e_{pq}+\lambda_{p,q \pm m} e_{p,q \pm m}) \in
\Jalg \langle H,x \rangle,
\]

\noindent where ${q\pm m}$ means $q+m$ if $q\in \{n+1, \dots , n+m\}$ and $q-m$ if
$q\in \{n+m+1, \dots , n+2m \}$ .

Assume $n>1$ and consider the element $(e_{qk} - q^{t}e_{kq}) \in
H(A,*)$ with $k \in \{1, \dots , n\}$ and $k \neq p$, then it
follows that $2(e_{qk} - q^{t}e_{kq}) \circ e_{pq} = e_{pk} \in
\Jalg \langle H,x \rangle$ with $p,k \in \{1, \dots ,n \}$ and $k
\neq p.$ Therefore we have found an element
 $\begin{pmatrix} a & 0 \\
0 & 0 \end{pmatrix} \in \Jalg \langle H,x \rangle $  with  $a \in
M_n(F)$ and $a \notin H(M_n(F), t)$ ($t$ denotes the usual transpose
involution). Since $H(M_n(F), t)$ is maximal subalgebra of
$M_n(F)^+$ (see  \cite[Theorem 6]{Ra1}) we obtain that
  $$\Jalg \langle H(M_n(F), t),a \rangle
 = M_n(F)^+$$

\noindent and so
$$\begin{pmatrix} M_n(F) & 0 \\
0 & 0 \end{pmatrix} \subseteq \Jalg \langle H,x \rangle.$$

\bigskip

Besides,  for any skewsymmetric matrix $ a \in M_n(F)$
 and for every $ b \in M_{n \times 2m}(F)$ one has
\[
\left[ \begin{pmatrix} a & 0 \\
0 & 0 \end{pmatrix} \circ \begin{pmatrix} 0 & b \\
-q^tb^t & 0 \end{pmatrix} \right] + \frac{1}{2} \begin{pmatrix} 0 & ab \\
-q^t(ab)^t & 0 \end{pmatrix} =  \begin{pmatrix} 0 & ab \\
0 & 0 \end{pmatrix} \in \Jalg \langle H,x \rangle,
\]

\noindent and thus
$\bigl(\begin{smallmatrix} 0 & M_{n \times 2m}(F) \\
0 & 0 \end{smallmatrix}\bigr) \subseteq \Jalg \langle H,x \rangle$,
because it is easy to check that
 $$K(M_n(F), t)M_{n \times
2m}(F)= M_{n \times 2m} (F).$$

\bigskip

But also
$$\begin{pmatrix} a & 0 \\ 0 & 0 \end{pmatrix} \circ
\begin{pmatrix} 0 & -b^tq^t \\ b & 0 \end{pmatrix} + \frac{1}{2}
\begin{pmatrix} 0 & -(ba)^tq^t \\ ba & 0 \end{pmatrix} \in
\Jalg \langle H,x \rangle
$$
and hence
\[
\begin{pmatrix} 0 & 0 \\
M_{2m \times n}(F) & 0 \end{pmatrix} \subseteq \Jalg \langle H,x
\rangle\qquad\text{and}\qquad
\begin{pmatrix} 0 & 0 \\
0 & M_{2m}(F) \end{pmatrix} \subseteq \Jalg \langle H,x \rangle.
\]

\bigskip

Finally, if $n=1$ then $\lambda e_{1j}+ \mu e_{1, j\pm m} \in \Jalg
\langle H, x \rangle $, with $j + m$ for $j \in \{n+1, \ldots ,
n+m\}$, and $j - m$ for $j \in \{n+m+1, \ldots , n+2m\}$. Now it is
clear that
\[
 \begin{pmatrix} M_n(F) & 0 \\ 0 & 0
\end{pmatrix}= \begin{pmatrix} F & 0 \\ 0 & 0
\end{pmatrix} \subseteq H(A,*) \subseteq \Jalg \langle H, x
\rangle.
\]

Taking $e_{j1} - e_{1, j \pm m} \in H$ one has
$$
2(\lambda e_{1j}+ \mu e_{1, j \pm m}) \circ (e_{j1}- e_{1, j \pm
m})= \lambda e_{11}+ \lambda e_{jj} \in \Jalg \langle H, x \rangle.
$$

\noindent Therefore, $e_{jj} \in \Jalg \langle H, x \rangle.$

Write $e_{jj}=\bigl(\begin{smallmatrix} 0 & 0 \\ 0 & a
\end{smallmatrix}\bigr)$ for a suitable $a \in M_{2m}(F)$. Then  $a \notin H(M_{2m}(F),
*)$ with $*$ the involution determined by the skewsymmetric bilinear
form with matrix
 $\bigl(\begin{smallmatrix} 0 & I \\ -I & 0
\end{smallmatrix}\bigr)$, and from the ungraded case (see \cite{Ra1})
we deduce that
 $$\Jalg \langle H(M_{2m}(F), *), a \rangle =
M_{2m}(F)^+$$

 \noindent and therefore $\begin{pmatrix} 0 & 0 \\ 0 & M_{2m}(F)
\end{pmatrix} \subseteq \Jalg \langle H, x \rangle$. Now it is easy to check
that since
$$\begin{pmatrix} 0 & b \\ -q^tb^t & 0
\end{pmatrix} \circ \begin{pmatrix} 0 & 0 \\ 0 & M_{2m}(F)
\end{pmatrix} \subseteq \Jalg \langle H, x \rangle$$
\noindent then $\begin{pmatrix} 0 & M_{1,2m}(F) \\
M_{1,2m}(F) & 0
\end{pmatrix} \subseteq \Jalg \langle H, x \rangle$ also in this case.

\medskip

If $x$ is now a nonzero homogeneous even element then
$$  x=\begin{pmatrix} a & 0 \\
0 & b \end{pmatrix}  $$

\noindent for a  skewsymmetric matrix $a$ and a matrix $b=
-q^tb^tq$. Consider
$$y=x \circ \begin{pmatrix} 0 & 0 \\
0 & I \end{pmatrix}=\begin{pmatrix} 0 & 0 \\
0 & b \end{pmatrix} \in \Jalg \langle H, x \rangle,$$

\noindent and
$$z= \begin{pmatrix} 0 & c \\
-q^tc^t & 0 \end{pmatrix}$$ \noindent such that $cb \neq 0.$ Then
$$ y  \circ  z = \frac{1}{2}\begin{pmatrix} 0 & cb \\
- bq^tc^t & 0 \end{pmatrix} \in \Jalg \langle H,x \rangle \cap
K_{\bar{1}}$$

\medskip

\noindent and the `odd' case applies.
\end{proof}

\bigskip

\subsection{$B = C^+$, $C\leq _{max}A$.}

Let us assume now that $B=C^+$ for a maximal subalgebra $C$  of the
simple associative superalgebra $A$. It has to be proved that $C^+$
is a maximal subalgebra of $A^+$.

Two different cases appear according to the classification of simple
associative superalgebras (see \cite {Wall}):

\begin{enumerate}

\item[\textbf{(1)}] $A$ is  simple as an (ungraded) algebra,
that is, $A$ is isomorphic to
$M_{p,q}(F)$, for some $p,q$. In this case, \cite[Theorem
2.2]{El-La-Sa} shows that either $C=eAe+ eAf+ fAf$ with $e,f$ even
orthogonal idempotents in $A$ such that $e+f=1$, or $C=C_A(u)$
(centralizer of $u$), with $u\in A_{\bar{1}}$ and $u^2=1$.

\item[\textbf{(2)}] $A$ is not  simple as an algebra, and hence it is
isomorphic to $Q_n(F)$ for some $n$. Then $A=A_{\bar{0}}+
A_{\bar{0}}u$ with $u\in Z(A)_{\bar{1}}$, $u^2=1$ and $A_{\bar{0}}$
is a simple algebra. In this case,  \cite[Theorem 2.5]{El-La-Sa}
shows that either $C=C_{\bar{0}}+ C_{\bar{0}}u$ with $C_{\bar{0}}$ a
maximal subalgebra of $A_{\bar{0}}$, or $C=A_{\bar{0}}$, or
$A_{\bar{0}}= D_{\bar{0}}+ D_{\bar{1}}$ is a $\bZ_2$-graded algebra
and $C=D_{\bar{0}}+ D_{\bar{1}}u$.

\end{enumerate}

\bigskip

\noindent\textbf{(1.a)}\quad Assume that $A$ is  simple as an
algebra, and that there are even orthogonal idempotents $e,f$ such
that $C=eAe+eAf+fAf$. Take an element $a_{\alpha}\in
A_{\alpha}\setminus C_{\alpha}$, so one has that $fa_\alpha e\neq
0$. Now  the element $(e \circ a_\alpha)\circ f = \frac{1}{4}
 (ea_\alpha f+fa_\alpha e)$ lies in
$\Jalg \langle C^+,a_\alpha \rangle$. Since $(fAf \circ fa_\alpha e)
\circ eAe = fAfa_\alpha eAe$, and $Afa_\alpha e A = A$, because $A$
is simple, it follows that $fAe\subseteq \Jalg \langle C^+,a_\alpha
\rangle$, and therefore $C^+$ is a maximal subalgebra of $A^+$. So
we have that in this case this condition is also sufficient to be a
maximal subalgebra of $A^+$.

\medskip

\noindent\textbf{(1.b)}\quad If $A$ is  simple as an algebra, but
$C= C_A(u)$, for an element $u\in A_{\bar{1}}$ with $u^2=1$, let $V$
be the irreducible $A$-module (unique, up to isomorphism), so that
$A$ can be identified with $\End_F(V)$. Then $u$ lies in
$\End(V)_{\bar{1}}$, and if $\{v_1, \dots, v_s\}$ is a basis of the
$F$-vector space $V_{\bar{1}}$, it follows that $\{ u(v_1), \dots,
u(v_s)\}$ is a $F$-basis of $V_{\bar{0}},$ and so $p=q$ and, since
$u^2=1$, the coordinate matrix of $u$ in this basis is
$$ u=\begin{pmatrix} 0 & I_s \\
I_s & 0 \end{pmatrix} .$$

\noindent Therefore $C_A(u) = Q_p(F)$, and then one can check easily
that $Q_p(F)$ is maximal in $M_{p,p}(F)$.

\medskip

\noindent\textbf{(2.a)}\quad Assume now that $A$ is not simple as an
algebra, so $A=A_{\bar{0}}+ A_{\bar{0}}u$, with $u\in
Z(A)_{\bar{1}}$, $u^2=1$ and $A_{\bar{0}}$  a simple algebra, and
that $C=C_{\bar{0}}+ C_{\bar{0}}u$, with $C_{\bar{0}}$ a maximal
subalgebra of $A_{\bar{0}}$. As for the ungraded case (see
\cite[page 192]{Ra1}) it follows that $\Jalg \langle
C_{\bar{0}}^+,a_{\bar{0}} \rangle=A_{\bar{0}}^+$ for any
$a_{\bar{0}}\in A_{\bar{0}} \setminus C_{\bar{0}}$. Thus
$A_{\bar{0}} \subseteq \Jalg \langle C^+,a_{\bar{0}} \rangle$.
Moreover since $1\in C_{\bar{0}},$ then $ u\in C$ and it follows
that $b_{\bar{0}}\circ u = \frac{1}{2} (b_{\bar{0}}u+ub_{\bar{0}})=
b_{\bar{0}}u\in \Jalg \langle C^+,a_{\bar{0}} \rangle$ for any
$b_{\bar{0}} \in A_{\bar{0}}$. Thus $A_{\bar{0}}u\subseteq \Jalg
\langle C^+,a_{\bar{0}} \rangle$ and $\Jalg \langle C^+,a_{\bar{0}}
\rangle=A^+$. Now take an element $a_{\bar{1}}\in A_{\bar{1}}
\setminus C_{\bar{1}}$. Then $a_{\bar{1}}=a_{\bar{0}}u$ with
$a_{\bar{0}}\in A_{\bar{0}} \setminus C_{\bar{0}}$. Since $u$ lies
in $C$, it follows that $a_{\bar{1}}\circ u = a_{\bar{0}} \in
 \Jalg  \langle C^+,a_{\bar{1}} \rangle $,
with $a_{\bar{0}}\in A_{\bar{0}}\setminus C_{\bar{0}}$ and the
`even' case applies.

\medskip

\noindent\textbf{(2.b)}\quad If $A$ is not simple as an algebra and
$C=A_{\bar{0}}$, let $b$ be any odd element: $b\in A\subuno =
A_{\bar{0}}u$. Thus $b=b_{\bar{0}}u$, for some $b_{\bar{0}}\in
A_{\bar{0}}$. Then $a_{\bar{0}}\circ b = (a_{\bar{0}}\circ
b_{\bar{0}})u$,  so $\Jideal \langle b_0 \rangle u\subseteq \Jalg
\langle A_{\bar{0}}^+,b \rangle$ (where $\Jideal\langle b_{\bar{0}}
\rangle$ denotes the ideal generated by $b_{\bar{0}}$ in the Jordan
algebra $A_{\bar{0}}^+$). By simplicity of $A_{\bar{0}}^+$,
$A_{\bar{0}}u \subseteq \Jalg \langle A_{\bar{0}}^+,b \rangle$, that
is, $C^+$ is a maximal subalgebra of $A^+$.

\medskip

\noindent\textbf{(2.c)}\quad Finally, assume that $A$ is not simple
as an algebra, and $A\subo$ (which is isomorphic to $M_p(F)$ for
some $p$) is $\bZ_2$-graded: $A\subo=D\subo\oplus D\subuno$, and
$C=D\subo\oplus D\subuno u$, where $u\in Z(A)\subuno$, $u^2=1$.
Here, as an associative superalgebra ($\bZ_2$-graded algebra),
$A\subo$ is isomorphic to $M_{r,s}(F)$ for some $r,s$. Identify
$A\subo$ to $M_{r,s}(F)$, so that $D\subo
=\left\{\bigl(\begin{smallmatrix}a&0
\\ 0&b\end{smallmatrix}\bigr): a\in M_r(F),\, b\in
M_s(F)\right\}$, and $D\subuno =\left\{\bigl(\begin{smallmatrix}0&u
\\ v&0\end{smallmatrix}\bigr): u\in M_{r\times s}(F),\, v\in
M_{s\times r}(F)\right\}$. Let us show that $C^+$ is a maximal
subalgebra of $A^+$. Since $A^+=C^+\oplus\bigl( D\subuno\oplus
D\subo u\bigr)$, it is enough to check that for any nonzero element
$x\in D\subo u\cup D\subuno$, the subalgebra of $A^+$ generated by
$C^+$ and $x$: $\Jalg\langle C^+,x\rangle$, is the whole $A^+$.

Take $0 \neq x \in D_{\bar{0}}u$. Then
$$ x=\begin{pmatrix} x_0 & 0 \\
0 & x_1 \end{pmatrix} u$$

\noindent with  $x_0 \in M_r(F),$ and  $x_1 \in M_s(F)$ not being
both zero. Without loss of generality, assume  $x_0 \neq 0$, and
take elements
$$ \begin{pmatrix} b & 0 \\
0 & 0 \end{pmatrix} \in C$$

\noindent with $0 \neq b\in M_r(F)$. Then
$$
\begin{pmatrix} b & 0 \\
0 & 0 \end{pmatrix}\circ x =  \begin{pmatrix} b & 0 \\
0 & 0 \end{pmatrix} \circ \begin{pmatrix} x_0 & 0 \\
0 & x_1 \end{pmatrix}u =\begin{pmatrix} b\circ x_0 & 0 \\
0 & 0 \end{pmatrix}u \in \Jalg \langle C^+,x \rangle $$ for any
$b\in M_r(F)$. Therefore
$$
\begin{pmatrix} \Jideal \langle x_0 \rangle & 0 \\
0 & 0 \end{pmatrix}u  \subseteq \Jalg  \langle C^+,x \rangle$$
\noindent and because of the simplicity of $M_n(F)^+$,
$$
\begin{pmatrix} M_r(F) & 0 \\
0 & 0 \end{pmatrix}u  \subseteq \Jalg  \langle C^+,x \rangle .$$

\medskip

Thus
\[  \begin{pmatrix} M_r(F) & 0 \\
0 & 0 \end{pmatrix}u \circ \begin{pmatrix} 0 & M_{r\times s}(F) \\
M_{s \times r}(F) & 0 \end{pmatrix}u \]
$$ = \begin{pmatrix} 0 & M_{r\times s}(F) \\
M_{s \times r}(F) & 0 \end{pmatrix} \subseteq \Jalg \langle C^+,x
\rangle,
$$
\noindent that is, $D_{\bar{1}} \subseteq \Jalg \langle C^+,x
\rangle$, and so $D_{\bar{1}} \circ D_{\bar{1}}u =
 D_{\bar{0}}u \subseteq \Jalg \langle C^+,x \rangle$ and $\Jalg  \langle C^+,x \rangle=A.$

\bigskip

Take now an element $0 \neq x\in D_{\bar{1}}$. Then an element
$d_{\bar{1}}u \in C^+$ can be found such that $0 \neq x \circ
d_{\bar{1}}u \in D_{\bar{0}}u \cap \Jalg \langle C^+,x \rangle$, so
the previous arguments apply.

\bigskip

This concludes the proof of the next result:

\begin{theo}\label{th:BCplus} Let $A$ be a finite dimensional simple
associative superalgebra over an algebraically closed field of
characteristic zero, and let $B$ be a maximal subalgebra of $A^+$
such that $B^\prime \neq A $ (where $B^\prime$ denotes the
associative subalgebra generated by $B$ in $A$). Then $B$ is a
maximal subalgebra of $A^+$ if and only if there is a maximal
subalgebra $C$ of the superalgebra $A$ such that $B= C^+$.
\end{theo}

\subsection{$B^\prime = A$ and $B$ is not semisimple.}

This situation does not appear in the ungraded case \cite{Ra1}.
However, consider the associative superalgebra $A=M_{1,1}(F)$ and
the subalgebra $B$ of $A^+$ spanned by
$\{e_{11},e_{22},e_{12}+e_{21}\}$, which, by dimension count, is
obviously maximal and satisfies that $B'=A$. The radical of $B$
consists of the scalar multiples of $e_{12}+e_{21}$, so it is
nonzero.

\medskip

\noindent\textbf{Question:} Is this, up to isomorphism, the only
possible example of a maximal subalgebra $B$ of $A^+$, $A$ being a
simple finite dimensional superalgebra over an algebraically field
$F$ of characteristic $0$, such that $B'=A$ and $B$ is not
semisimple?

\bigskip

\section {Maximal subalgebras of $H(A,*).$}

\bigskip

Consider now the Jordan superalgebra $J=H(A,*)$, where $A$ is a
finite dimensional simple associative superalgebra over an
algebraically closed field $F$ of characteristic zero, and $*$ is a
superinvolution of $A$.

Up to isomorphism  \cite[Theorem 3.1]{Go-She}, it is known that
$A=M_{p,q}(F)$ and that $*$ is either the orthosymplectic
 or the transpose superinvolution, that is, $H(A,*)$
is either  $osp_{n,2m}$ or $p(n)$.

Let $B$ be a maximal subalgebra of $H(A,*),$ then again three
possible situations appear:

\begin{enumerate}

\item[{\rm (i)}] either $B^\prime = A$  and $B$ is semisimple,

\item[{\rm (ii)}] or $B^\prime \neq A,$

\item[{\rm (iii)}] or $B^\prime =A$ and $B$ is not semisimple.

\end{enumerate}

\subsection{$B^\prime = A$ and $B$ semisimple.}

Let us assume first that $B$ is a maximal subalgebra of the simple
superalgebra $H(A,*)$, with $B'=A$ and $B$ semisimple. From Lemma
\ref{le:BinAplusimproved}, we know  that either $B$ is isomorphic to
$D_t$ ($t\ne 0,\pm 1,-2,-\frac{1}{2}$), or $B=H(A,\diamond)$ with
$\diamond$ a superinvolution. In the first case we remark that we
have given only necessary conditions in Proposition 3.3 if $B^\prime
=A$ and $1_A\in B.$ In the second case, one has $B=H(A,\diamond)
\subseteq H(A,*),$, but  Theorem \ref{th:BsemiAplus} shows that
$H(A,\diamond)$ is maximal in $A^+$, thus obtaining a contradiction.

\bigskip

Therefore:

\begin{theo}\label{th:BsemiinHAstar}
Let $J$ be the Jordan superalgebra  $H(A, *)$, where $A$ is a finite
dimensional simple Jordan superalgebra over an algebraically closed
field of characteristic zero, and $*$ a superinvolution in $A$. If
$B$ is a maximal subalgebra of $J$ such that $B^{\prime}=A$ and $B$
is semisimple, then $B= D_t$ ($t\ne 0,\pm 1,-2,-\frac{1}{2}$) and
$(A,*)$ is given by  Proposition \ref{pr:DtEndV}.
\end{theo}

\medskip

\noindent\textbf{Question:} Given a natural number $m$, and with $t$
equal either to $-\frac{m}{m+1}$ or to $-\frac{m+1}{m}$, is $D_t$
isomorphic to a maximal subalgebra of the Jordan superalgebra
$H(\End_F(V),*)$ ($V$ and $*$ as in Proposition \ref{pr:DtEndV})?

\medskip

For $m=2$ or $m=3$, this has been checked to be the case.

\bigskip

\subsection{$B'\ne A$.}

Assume now that the maximal subalgebra $B$ of $H(A,*)$ satisfies
$B'\ne A$. The result that settles this case is the following:

\bigskip

\begin{theo}\label{th:BprimeneAinHAstar}
Let $J$ be the Jordan superalgebra $H(A,*)$, where $A$ is a finite
dimensional simple Jordan superalgebra over an algebraically closed
field of characteristic zero, and $*$ is a superinvolution in $A$.
Let $B$ be a subalgebra of $J$ such that $B^\prime \neq A$ (where as
always
 $B^\prime $ is the subalgebra of $A$ generated by $B$). Then $B$ is maximal if and only
if there are even idempotents $e,f\in A$ with $e+f=1$ such that
$B=H(C,*)$ and one of the following possibilities occurs:

\begin{enumerate}

\item[{\rm (i)}] either  $C=eAe+fAf$, $e^*=e$, $f^*=f$,  $H(eAe,*)^\prime
 =eAe$, and $H(fAf,*)^\prime =fAf$.

\item[{\rm (ii)}] or $C=eA +Ae^* + ff^*Aff^*$, with $H(ff^*Aff^*,*)^\prime =
ff^*Aff^*$.

\end{enumerate}

\end{theo}

Note \cite{Go} that given a finite dimensional simple associative
superalgebra $C$ over $F$ with a superinvolution $*$, the
associative subalgebra $H(C,*)'$ is the whole $C$ unless $(C,*)$ is
either a quaternion superalgebra with the transpose superinvolution
or a quaternion algebra with the standard involution.

\begin{proof} If $B^\prime = A$, and since $B\subseteq H(A,*)$, it follows that
$B^\prime$ is closed under the superinvolution $*$, and so $B^\prime
\subseteq C$ with $C$ a maximal subalgebra of $(A,*)$. But using the
maximality of $B$ and that $B\subseteq H(A,*)$, one concludes  that
 $B=H(C,*)$. Recall that  $H(A,*)$ is isomorphic either to $p(n)$ or to
 $osp_{n,2m}$.

If $B=H(C,*)$ with $C$ a maximal subalgebra of $(A,*)$, then the
results in \cite{El-La-Sa} show that either $C= (eAe+eAf+fAf)\cap
(e^*Ae^*+f^*Ae^*+f^*Af^*)$ with $e,f$ even orthogonal idempotents,
or $C= C_A(u)$ with $u\in A_{\bar{1}},$ $0\ne u^2\in F,$ $u^*\in
Fu$. In this last case, since $u^* \in Fu$ it follows that $u^* =
\alpha u$ with $\alpha \in F$. But $(u^*)^* =u$ and so $\alpha^2
=1$, that is, $\alpha = \pm 1$. Thus $u^2=(u^2)^*=-(u^*)^2= -u^2,$ a
contradiction.

Thus,  $C$ is of the first type, and then  \cite[Proposition
4.6]{El-La-Sa} gives two possible cases.

In the first case there is an idempotent $e$ of $A$ such that
$C=eAe+fAf$ and $e^*=e,$ $f=1-e.$ If $H(C,*)^\prime \neq C$ then
either $H(eAe,*)^\prime \neq eAe$ or $H(fAf,*)^\prime\neq fAf$. It
may be assumed that  $H(eAe,*)^\prime \neq eAe$, and then the
results in \cite{Go} show that either $eAe$ is a quaternion
superalgebra with the restriction $*|_{eAe}$ being the transpose
superinvolution or is a quaternion algebra contained in
$A_{\bar{0}}$, with the standard involution. In both cases
$e=e_1+e_2$ with $e_1,e_2$ orthogonal idempotents and $e_1^*=e_2.$
Consider $D=e_1A+Ae_2+fAf$ and take  $0\neq e_1af \in e_1Af$, then
$e_1af+fa^*e_2\in H(D,*)$ and $e_1af+fa^*e_2\notin H(C,*)$. In the
same vein, take $c\in A$ with $e_2cf\neq 0$. Then $e_2cf+fc^*e_1 \in
H(A,*)\setminus H(D,*).$ Therefore $B=H(C,*) \subsetneqq H(D, *)
\subsetneqq H(A,*)$ and $B=H(C,*)$ is not maximal. So $B^\prime=
H(C,*)^\prime =C$ if $B=H(C,*)$ with $C=eAe+ fAf$ and $e^*=e.$

In the second case  \cite[Proposition 4.6]{El-La-Sa}, there is an
idempotent $e$ in $A$ such that $e$, $e^*$, $ff^*$ are mutually
orthogonal idempotents with $1=e+e^*+ff^*$, and $C=
eA+Ae^*+ff^*Aff^*$. Hence $H(C,*) = H(ff^*aff^*) + \{ ea + a^*e^* :
a\in A\}$.

If $H(ff^*Aff^*,*)^\prime \neq ff^*Aff^*$, then $ff^*Aff^*$ is a
quaternion superalgebra  with superinvolution such that $ff^* =
e_1+e_2$ with $e_1,e_2$ orthogonal idempotents and $e_1^*=e_2$.
Consider the subalgebra $D = eA+Ae^*+e_2A+Ae_1$. As
$H(C,*)\subsetneqq H(D,*)\subsetneqq H(A,*)$,  $H(C,*)$ is not
maximal. Therefore, if $B=H(C,*)$ with $C=eA+Ae^*+ff^*Aff^*$, and
$e$, $e^*$, $ff^*$ mutually orthogonal idempotents such that
$e+e^*+ff^*=1$, then $H(ff^*Aff^*, *)^\prime =ff^*Aff^*$.

\bigskip

The proof of the converse will be split according to the different
possibilities:

\medskip

\noindent\textbf{(i.1)}: The superinvolution $*$ on $A$ is the
transpose superinvolution, and the conditions in item (i) of the
Theorem hold:

\smallskip

Then $*$ is determined, after identifying  $A$ with $\End_F(V)$, by
a nondegenerate odd symmetric superform $( \ , \ )$. That is, ,
$(V_{\bar{0}}, V_{\bar{0}}) = (V_{\bar{1}}, V_{\bar{1}})=0$ and
$(a_0, b_1)= (b_1, a_0)$ for any $a_0 \in V_{\bar{0}}$, $b_1 \in
V_{\bar{1}}$.

In this situation we claim that  a basis $\{x_1, \ldots ,x_n,y_1
\ldots , y_n\}$ of $ V$ can be chosen such that $\{x_1, \ldots
,x_n\}$ is a basis of $V_{\bar{0}}$, $\{y_1, \ldots ,y_n\}$ is a
basis of $V_{\bar{1}}$, and the coordinate matrices of the superform
and of $e$ present the following form, respectively,
$$\begin{pmatrix} 0 & 0 & I &  0  \\
0 & 0 & 0 &  I  \\ I & 0 & 0 &  0  \\ 0 & I & 0 & 0 \end{pmatrix},
\hspace{1cm} \begin{pmatrix} I & 0 & 0 & 0  \\ 0 & 0 & 0 & 0 \\ 0 &
0 & I & 0 \\ 0 & 0 & 0 & 0 \end{pmatrix}.$$

\noindent This follows from the fact that the eigenspaces of the
idempotent transformation $e$ are orthogonal relative to $(\ ,\ )$,
as $e^*=e$. Under these circumstances, we may identify $H(A,*)$ to
\[
p(n)=\left\{\begin{pmatrix} a&b\\ c&a^t\end{pmatrix} : \text{$b$
skewsymmetric, $c$ symmetric}\right\}
\]
in such a way that the subalgebra $H(eAe+fAf,*)$ becomes the
subspace of the matrices (in block form)
\[
\begin{pmatrix}
a_{1}&0&c_{1}&0\\ 0&a_{2}&0&c_{2}\\ d_{1}&0&a_{1}^t&0\\
 0&d_{2}&0&a_{2}^t\end{pmatrix}
\]
where $a_{1},c_{1},d_{1}$ belong to $M_i(F)$, $a_{2},c_{2},d_{2}$
belong to $M_j(F)$, $i+j=n$, and $c_{1},c_2$ are skewsymmetric
matrices, while $d_1,d_2$ are symmetric.

It must be proved that for any homogeneous element $x$, $\Jalg
\langle H(C,*),x \rangle =H(A,*)$ holds.

\bigskip

Let $x \in H(A,*)_{\bar{0}} \setminus H(C,*)_{\bar{0}}$, that is,
\[
x=\sum_{\substack{1\leq k \leq i \\[0.1cm] 1 \leq r \leq j}}
 \lambda_{kr} (e_{k,i+r}+e_{n+i+r,n+k})
 + \sum_{\substack{1\leq r \leq j \\[0.1cm] 1 \leq k \leq i}} \mu_{rk} (e_{i+r,k}+e_{n+k,n+i+r})
\]
where $e_{r,s}$ denotes the matrix with $1$ in the $(r,s)$-th entry
and  $0$
 in all the other entries. Suppose that there exists $\lambda_{pq} \neq 0.$
The same proof works if $\mu_{pq} \neq 0$.

Since $H(C,*)'=C$ and $i>1$ (as $H(eAe,*)'=eAe$), an index  $s \in
\{1, \ldots, i\}$ can be chosen with $s \neq p$, such that
$u=e_{s,p}+e_{n+p,n+s} \in H(C,*)$. Let $v=e_{p,p} + e_{n+p,n+p}$
and $w= e_{i+q,i+q} + e_{n+i+q,n+i+q}$ (note that
 $v,w \in H(C,*)$). Then
$$((v \circ x) \circ w )\circ u= \frac{1}{8} \lambda_{pq}(e_{s,i+q}+e_{n+i+q,n+s})
\in \Jalg \langle H(C,*),x \rangle.$$

Denote this element by $\alpha$, and then $0 \neq \alpha \in
e_1Af_1+ f_1^*Ae_1^*.$ Now

$$((e_1ae_1+e_1^*a^*e_1^*) \circ \alpha ) \circ (f_1bf_1+f_1^*b^*f_1^*)=
 e_1ae_1\alpha f_1bf_1
+f_1^*b^*f_1^*\alpha e_1^*a^*e_1^* $$

\noindent belongs to $ \Jalg \langle H(C,*),x \rangle$. Since
$\{ae_1 \alpha f_1 b :   a,b \in A \} $ is an ideal of $ A$, and $A$
is simple, it holds that
 $\{ae_1 \alpha f_1 b :
a,b \in A\}=A$, and so $e_1 a f_1 + f_1^* a^* e_1^* \in \Jalg
 \langle H(C,*),x \rangle$ for any $a \in A.$

Consider now an element $y \in f_1Af_1^* \cap H(C,*)$. Since $j>1$
(because $ H(fAf,*)'=fAf$), we can pick up the element $ y=
e_{l,k}-e_{l+1,k-1}$, with $l=i+1$ and $k=n+i+2$. Take
$z=e_{k-1,p}+e_{1,l} \in H(e_1Af_1 + f_1^* A e_1^*,*) \subseteq
\Jalg \langle H(C,*),x \rangle$ and $v= e_{p,1} \in H(C,*) \cap
e_1^*Ae_1$, with $p=n+1.$ Then $(y \circ z) \circ v=
\frac{1}{4}(-e_{l+1,1}-e_{p,k}) \in (f_1Ae_1 +e_1^*Af_1^*) \cap
H(A,*)_{\bar{0}}$. As before we obtain that $f_1ae_1 +e_1^*a^* f_1^*
\in \Jalg \langle H(C,*),x \rangle$, and $H(A,*)_{\bar{0}} \subseteq
\Jalg \langle H(C,*),x \rangle$.

Now it will proved that $H(A,*)_{\bar{1}}$ is contained in $\Jalg
\langle H(C,*), x \rangle $. Take $y=e_{k, n+i+t}- e_{i+t,n+k} \in
H(A,*)_{\bar{1}} \cap (e_1Af_1^*+f_1Ae_1^*)$, with $k \in \{1,
\ldots ,i\}, t \in \{1, \ldots ,j\}$ and we claim that $y \in \Jalg
\langle H(C,*),x \rangle$. Since $H(fAf,*)^{\prime}= fAf,$ there
exists $s \in \{1, \ldots ,j\}$ with $s \neq t$,
 and consider then the elements
$z=e_{n+i+s,n+k}+ e_{k,i+s} \in \Jalg \langle H(C,*),x \rangle$, and
$u=e_{i+s,n+i+t}-e_{i+t,n+i+s} \in H(C,*)$. Then it follows that $z
\circ u=\frac{1}{2}y \in \Jalg \langle H(C,*),x \rangle$. In the
same way we obtain that $(e_1^*Af_1 + f_1^*Ae_1) \cap
H(A,*)_{\bar{1}} \subseteq \Jalg \langle H(C,*),x \rangle$.

So for any $x \in H(A,*)_{\bar{0}}\setminus H(C,*)_{\bar{0}}$,
$H(A,*)=\Jalg \langle H(C,*), x \rangle$ holds.

\bigskip

Now let $x \in H(A,*)_{\bar{1}} \setminus H(C,*)_{\bar{1}}.$ Then
$$x= \sum_{\substack{1\leq k \leq i \\[0.1cm] 1 \leq r \leq j}}
\lambda_{kr} (e_{k,n+i+r}-e_{i+r,n+k})+
\sum_{\substack{1\leq k \leq i \\[0.1cm] 1 \leq r \leq j}}
\mu_{kr}(e_{n+k,i+r}+e_{n+i+r,k})$$

\noindent and assume that  for some $(p,q)$, one has $\lambda_{pq}
\neq 0$.

Since $u=e_{n+p,p} \in H(C,*),$ $0 \neq 2 x \circ u= - \sum_{1\leq q
\leq j} \lambda_{pq} (e_{i+q,p}+e_{n+p,n+i+q})  \in
H(A,*)_{\bar{0}}\setminus H(C,*)_{\bar{0}}$, and  the above case
applies.

In the same way, if $\mu_{pq} \neq 0$ we obtain that $H(C,*)$ is a
maximal subalgebra of $H(A,*).$

\bigskip

\noindent\textbf{(i.2)}: The superinvolution $*$ on $A$ is an
orthosymplectic superinvolution, and the conditions in item (i) of
the Theorem hold:

\smallskip

In this and the following cases, we will content ourselves to
establish the setting in which one can apply the same kind of not
very illuminating arguments like those used in case \textbf{(i.1)}.

Here, after identifying $A$ to $\End_F(V)$, the superinvolution $*$
is determined by a nondegenerate symmetric superform $(\ , \ )$ on
$V$, that is,
 $(\ , \ )\vert_{V_{\bar{0}} \times V_{\bar{0}}}$ is symmetric,
$(\ , \ )\vert_{V_{\bar{1}} \times V_{\bar{1}}}$ is skewsymmetric
and $(V_{\bar{0}},V_{\bar{1}})=(V_{\bar{1}}, V_{\bar{0}})=0$.

Since $e$ is idempotent and selfadjoint, there is a basis of $V$ in
which the coordinate matrices of the superform and of $e$ are,
respectively,
\[
\begin{pmatrix} I&0&0&0&0&0\\
   0&I&0&0&0&0\\
   0&0&0&0&I&0\\
   0&0&0&0&0&I\\
   0&0&-I&0&0&0\\
   0&0&0&-I&0&0\end{pmatrix},\qquad
\begin{pmatrix} I&0&0&0&0&0\\
                0&0&0&0&0&0\\
                0&0&I&0&0&0\\
                0&0&0&0&0&0\\
                0&0&0&0&I&0\\
                0&0&0&0&0&0
                \end{pmatrix},
\]
where $0$, respectively $I$, denotes the zero matrix, respectively
identity matrix (of possibly different orders). Let $n$ be the
dimension of $V\subo$, $2m$ the dimension of $V\subuno$, $i$ the
rank of the restriction $e\vert_{V\subo}$, $j=n-i$,  $2k$ the rank
of $e\vert_{V\subuno}$ and $l=m-k$. Hence, identifying by means of
this basis $H(A,*)$ to $osp_{n,2m}$, the idempotent $e$ decomposes
as $e= e_1+e_2+ e_2^*$, with $e_1= \sum_{s=1}^i e_{s,s}$,
$e_2=\sum_{s=1}^k e_{n+s,n+s}$ and $e_2^*=\sum_{s=1}^k
e_{n+m+s,n+m+s}$. Similarly, $f=1-e$ decomposes as $f=f_1+ f_2+f_2^*
$.

The elements of
 $H(C,*)$ are then the matrices (in block form)
\[
\left( \begin{array} {cccccc} c_{11} & 0 \hspace{0.8cm} \vdots & b_{11} & 0 & b_{13} & 0\\
0 & c_{22} \  \hspace{0.25cm}  \vdots  & 0 & b_{22} & 0 & b_{24}\\
\hdotsfor{6} \\
b_{13}^t & 0 \hspace{0.80cm} \vdots &a_{11} & 0 & a_{13}  & 0 \\
0 & b_{24}^t \  \hspace{0.25cm} \vdots &0 & a_{22} & 0  & a_{24} \\
-b_{11}^t & 0 \hspace{0.80cm} \vdots &a_{31} & 0 & a_{11}^t  & 0 \\
0 & -b_{22}^t \  \vdots & 0 & a_{42} & 0  & a_{22}^t
\end{array}\right)
\]
 with $c_{11} \in M_i(F) $ and $ c_{22} \in M_j(F)$ symmetric matrices,
$a_{11} \in M_k(F)$, $a_{22} \in M_l(F)$, $b_{11}, b_{13} \in
M_{i\times k}(F)$, $b_{22}, b_{24} \in M_{j\times l}(F),$ $a_{13},
a_{31} \in M_k(F)$ skewsymmetric matrices, and $a_{24}, a_{42} \in
M_l(F) $  skewsymmetric too

Note that it is possible that either $e_1$ or $f_1$ may be $0$. If,
for instance, $f_1=0$, then since $H(fAf,*)'=fAf$, it follows that
$l>1$.

In this setting, routine arguments like the ones for \textbf{(i.1)}
apply.

\bigskip

\noindent\textbf{(ii.1)}: The superinvolution $*$ on $A$ is the
transpose superinvolution, and the conditions in item (ii) of the
Theorem hold:

\smallskip

Here a basis $\{x_1, \ldots ,x_n,y_1 \ldots , y_n\}$ of $V$ ($\{x_1,
\ldots ,x_n\}$ being a basis of $V_{\bar{0}}$ and $\{y_1, \ldots
,y_n\}$  of $V_{\bar{1}}$), so that the coordinate matrices of the
superform and of the idempotents $e$, $e^*$ and $ff^*$ are,
respectively,

$$\left(\begin{array} {cccccc} 0 & 0 & 0 & I &  0  &  0\\
0 & 0 & 0 & 0 &  I  & 0 \\
0 & 0 & 0 & 0 &  0  & I \\
I & 0 & 0 &  0  &  0 &  0\\
0 & I & 0 & 0 & 0 & 0 \\
0 & 0 & I & 0 & 0 & 0\end{array}\right),
 \hspace{1cm} \left( \begin{array} {cccccc} I & 0 & 0 & 0  & 0 & 0 \\
0 & 0 & 0 & 0 & 0 & 0  \\
0 & 0 & 0 & 0 & 0 & 0  \\
0 & 0 & 0 & 0 & 0 & 0  \\
0 & 0 & 0 & 0 & 0 & 0  \\
0 & 0 & 0 & 0 & 0 & I\end{array} \right),$$

$$\left( \begin{array} {cccccc} 0 & 0 & 0 & 0  & 0 & 0 \\
0 & 0 & 0 & 0 & 0 & 0  \\
0 & 0 & I & 0 & 0 & 0  \\
0 & 0 & 0 & I & 0 & 0  \\
0 & 0 & 0 & 0 & 0 & 0  \\
0 & 0 & 0 & 0 & 0 & 0\end{array} \right), \hspace{1cm}
\left( \begin{array} {cccccc}0 & 0 & 0 & 0  & 0 & 0 \\
0 & I & 0 & 0 & 0 & 0  \\
0 & 0 & 0 & 0 & 0 & 0  \\
0 & 0 & 0 & 0& 0 & 0  \\
0 & 0 & 0 & 0 & I & 0  \\
0 & 0 & 0 & 0 & 0 & 0\end{array} \right).$$

This follows from the fact that $e$, $e^*$ and $ff^*$ are orthogonal
idempotents with $1=e+e^*+ff^*$, so
\[
\begin{split}
V\subo&=S(1,e)\subo\oplus S(1,ff^*)\subo\oplus S(1,e^*)\subo,\\
V\subuno&=S(1,e^*)\subuno\oplus S(1,ff^*)\subuno\oplus
S(1,e)\subuno,
\end{split}
\]
where $S(1,g)$ denotes the eigenspace of the endomorphism $g$ of
eigenvalue $1$, and from the fact that  $ff^*$ is selfadjoint, so
\[
V=\bigl(S(1,e)\subo\oplus S(1,e^*)\subuno\bigr)\oplus
 \bigl(S(1,ff^*)\subo\oplus S(1,ff^*)\subuno\bigr)\oplus
 \bigl(S(1,e^*)\subo\oplus S(1,e)\subuno\bigr).
\]

After the natural identifications, the elements of
$H(C,*)=H(eA+Ae^*+ff^*Aff^*,*)$ are the matrices (in block form)
\[
\begin{pmatrix} a_{11} & a_{12}  & a_{13} \thickspace \vdots& c_{11} & c_{12} & c_{13} \\
  0 & a_{22}  & a_{23} \thickspace \vdots & -c_{12}^t & c_{22}& 0 \\
  0 & 0       & a_{33} \thickspace \vdots& -c_{13}^t & 0     & 0 \\
\hdotsfor{6}\\
  0 & 0       & d_{13} \thickspace \vdots& a_{11}^t  & 0     & 0 \\
  0 & d_{22}  & d_{23} \thickspace \vdots& a_{12}^t  & a_{22}^t & 0 \\
  d_{13}^t & d_{23}^t  & d_{33}  \thickspace \vdots & a_{13}^t  &
      a_{23}^t & a_{33}^t  \end{pmatrix}
\]
where $c_{11}$, $c_{22}$ are skewsymmetric matrices and $d_{22}$,
$d_{33}$ symmetric matrices. Since
$H(ff^*Aff^*,*)^{\prime}=ff^*Aff^*$, it follows that $ff^*Aff^*$ is
not a quaternion superalgebra and so the order of the blocks in the
$(2,2)$ position is $>1$.

This is the setting where routine computations can be applied.

\bigskip

\noindent\textbf{(ii.2)}: The superinvolution $*$ on $A$ is an
orthosymplectic superinvolution, and the conditions in item (ii) of
the Theorem hold:

\smallskip

Here, with the same sort of arguments as before, the coordinate
matrices in a suitable basis of the orthosymplectic superform, and
of the idempotents $ff^*$, $e$ and $e^*$ are, respectively:
\[
\left(\begin{smallmatrix}
I & 0 & 0 & 0 &  0 & 0 & 0\\
0 & 0 & I & 0 & 0 &  0  & 0 \\
0 & I & 0 & 0 & 0 &  0  & 0 \\
0 & 0 & 0 & 0 & 0 &  I  & 0 \\
0 & 0 & 0 & 0 & 0 &  0  & I \\
0 & 0 & 0 & -I &  0&  0  & 0 \\
0 & 0 & 0 & 0 & -I & 0 & 0 \end{smallmatrix}\right),
\hspace{0.5cm} \left(\begin{smallmatrix}
I & 0 & 0 & 0 &  0 & 0 & 0\\
0 & 0 & 0 & 0 & 0 &  0  & 0 \\
0 & 0 & 0 & 0 & 0 &  0  & 0 \\
0 & 0 & 0 & I & 0 &  0  & 0 \\
0 & 0 & 0 & 0 & 0 &  0  & 0 \\
0 & 0 & 0 & 0 & 0 &  I  & 0 \\
0 & 0 & 0 & 0 & 0 & 0 & 0 \end{smallmatrix}\right), \hspace{0.5cm}
\left( \begin{smallmatrix}
0 & 0 & 0 & 0 & 0 & 0 & 0\\
0 & I & 0 & 0 & 0 & 0 & 0 \\
0 & 0 & 0 & 0 & 0 &  0  & 0 \\
0 & 0 & 0 & 0 & 0 &  0  & 0 \\
0 & 0 & 0 & 0 & I &  0  & 0 \\
0 & 0 & 0 & 0 & 0 &  0  & 0 \\
0 & 0 & 0 & 0 & 0 & 0 & 0\end{smallmatrix} \right), \hspace{0.5cm}
\left( \begin{smallmatrix}
0 & 0 & 0 & 0 & 0 & 0 & 0\\
0 & 0 & 0 & 0 & 0 & 0 & 0 \\
0 & 0 & I & 0 & 0 &  0  & 0 \\
0 & 0 & 0 & 0 & 0 &  0  & 0 \\
0 & 0 & 0 & 0 & 0 &  0  & 0 \\
0 & 0 & 0 & 0 & 0 &  0  & 0 \\
0 & 0 & 0 & 0 & 0 & 0 & I\end{smallmatrix} \right).
\]

Now, the superinvolution $*$, identifying the elements in $H(A,*)$
with their coordinate matrices in the basis above, is given by:
\[
\left( \begin{smallmatrix} a_{11} & a_{12} & a_{13} & a_{14} &
a_{15}  & a_{16} & a_{17}
\vspace{0.15cm}  \\
a_{21} & a_{22} & a_{23} & a_{24} & a_{25}  & a_{26} & a_{27} \vspace{0.15cm} \\
a_{31} & a_{32} & a_{33} & a_{34} & a_{35}  & a_{36} & a_{37} \vspace{0.15cm} \\
a_{41} & a_{42} & a_{43} & a_{44} & a_{45}  & a_{46} & a_{47} \vspace{0.15cm} \\
a_{51} & a_{52} & a_{53} & a_{54} & a_{55}  & a_{56} & a_{57} \vspace{0.15cm} \\
a_{61} & a_{62} & a_{63} & a_{64} & a_{65}  & a_{66} & a_{67} \vspace{0.15cm} \\
a_{71} & a_{72} & a_{73} & a_{74} & a_{75}  & a_{76} & a_{77}   \end{smallmatrix} \right) \rightarrow \left( \begin{smallmatrix}
a_{11}^t & a_{31}^t & a_{21}^t & a_{61}^t & a_{71}^t  & -a_{41}^t & -a_{51}^t \\
a_{13}^t & a_{33}^t & a_{23}^t  & a_{63}^t & a_{73}^t  & -a_{43}^t & -a_{53}^t \\
a_{12}^t & a_{32}^t & a_{22}^t & a_{62}^t & a_{72}^t  & -a_{42}^t & -a_{52}^t \\
-a_{16}^t & -a_{36}^t & -a_{26}^t & a_{66}^t & a_{76}^t  & -a_{46}^t & -a_{56}^t \\
-a_{17}^t & -a_{37}^t & -a_{27}^t & a_{67}^t & a_{77}^t  & -a_{47}^t & -a_{57}^t \\
a_{14}^t & a_{34}^t & a_{24}^t & -a_{64}^t & -a_{74}^t  & a_{44}^t & a_{54}^t \\
a_{15}^t & a_{35}^t & a_{25}^t & -a_{65}^t & -a_{75}^t  & a_{45}^t &
a_{55}^t   \end{smallmatrix} \right).
\]
 Therefore the Jordan superalgebra $H(A,*)$ consists of the
following matrices:
\[
\begin{pmatrix}
a_{11} & a_{12} & a_{13} \hspace{0.4cm} \vdots  & a_{14} & a_{15}  & a_{16} & a_{17} \\
a_{13}^t & a_{22} & a_{23} \hspace{0.4cm} \vdots  & a_{24} & a_{25}  & a_{26} & a_{27} \\
a_{12}^t & a_{32} & a_{22}^t \hspace{0.4cm} \vdots & a_{34} & a_{35}  & a_{36} & a_{37} \\
\hdotsfor{7}\\
-a_{16}^t & -a_{36}^t & -a_{26}^t\vdots  & a_{44} & a_{45}  & a_{46} & a_{47} \\
-a_{17}^t & -a_{37}^t & -a_{27}^t \vdots & a_{54} & a_{55}  & -a_{47}^t & a_{57} \\
a_{14}^t & a_{34}^t & a_{24}^t \hspace{0.4cm} \vdots & a_{64} & a_{65}  & a_{44}^t & a_{54}^t \\
a_{15}^t & a_{35}^t & a_{25}^t \hspace{0.4cm} \vdots & -a_{65}^t &
a_{75}  & a_{45}^t & a_{55}^t  \end{pmatrix},
\]
where $a_{11}, a_{23}, a_{32} $ are symmetric matrices, while
$a_{46}, a_{57}, a_{64}, a_{75}$ are skewsymmetric matrices.
Besides, the elements of $H(C,*)=H(eA+Ae^*+ff^*Aff^*,*)$ are the
matrices which, in block form, look like
\[
\begin{pmatrix}
a_{11} & 0 & a_{13} \hspace{0.4cm} \vdots  & a_{14} & 0  & a_{16} & a_{17} \\
a_{13}^t & a_{22} & a_{23} \hspace{0.4cm} \vdots  & a_{24} & a_{25}  & a_{26} & a_{27} \\
0 & 0 & a_{22}^t \hspace{0.4cm} \vdots & 0 & 0  & 0 & a_{37} \\
\hdotsfor{7}\\
-a_{16}^t & 0 & -a_{26}^t\vdots  & a_{44} & 0  & a_{46} & a_{47} \\
-a_{17}^t & -a_{37}^t & -a_{27}^t \vdots & a_{54} & a_{55}  & -a_{47}^t & a_{57} \\
a_{14}^t & 0 & a_{24}^t \hspace{0.4cm} \vdots & a_{64} & 0  & a_{44}^t & a_{54}^t \\
0 & 0 & a_{25}^t \hspace{0.4cm} \vdots & 0 & 0  & 0 & a_{55}^t  \end{pmatrix}.
\]
Now again routine arguments with matrices give the result.
\end{proof}

\bigskip

\subsection{$B^\prime = A$ and $B$ is not semisimple.}

As for the maximal subalgebras of the Jordan superalgebras $A^+$,
this situation does not appear in the ungraded case \cite{Ra1}.
However, consider the associative superalgebra $A=M_{1,2}(F)$, with
the natural orthosymplectic superinvolution. Thus, the Jordan
superalgebra $J=H(A,*)$ is
\[
J=osp_{1,2}=\left\{\begin{pmatrix} a&-c&b\\ b&d&0\\
c&0&d\end{pmatrix} : a,b,c,d\in F\right\}.
\]
The subspace
\[
B=\left\{\begin{pmatrix} a&-b&b\\ b&d&0\\ b&0&d\end{pmatrix}:
a,b,d\in F\right\}
\]
is a maximal superalgebra of $J$, and it satisfies $B'=A$, while it
is not semisimple, as its radical coincides with its odd part

\medskip

\noindent\textbf{Question:} Is  this, up to isomorphism, the only
possible example of a maximal subalgebra $B$ of $H(A,*)$, $A$ being
a simple finite dimensional superalgebra over an algebraically field
$F$ of characteristic $0$, such that $B'=A$ and $B$ is not
semisimple?

\bigskip

It seems that a broader knowledge of non semisimple Jordan
superalgebras is needed here.

The solution to the above question is also related to the Question
after Theorem \ref{th:BsemiinHAstar}. Actually, if this question is
answered in the affirmative, then the subalgebra $B$ isomorphic to
$D_t$ ($t\ne 0,\pm 1,-2,-\frac{1}{2}$) in Theorem
\ref{th:BsemiinHAstar} would indeed be maximal in $H(A,*)$.
Otherwise, any maximal subalgebra $S$ containing $B$ would satisfy
$S'=A$ (as $B'=A$ already) and would not be semisimple (because of
Theorem \ref{th:BsemiinHAstar}).

\end{document}